%% file: hypergeom.tex
\documentclass{article}

\usepackage{amsfonts}
\usepackage{epsfig}
\usepackage{amsmath}
\usepackage{amssymb}
\usepackage{amsthm}
\usepackage{oldgerm}
\usepackage{subfigure}
\newcommand{\oz}{\langle p,q\rangle_{10}}
\newcommand{\zinf}{\langle p,q\rangle_{0\infty}}
\newcommand{\oinf}{\langle p,q\rangle_{1\infty}}
\newcommand{\F}{{_2F_1}}
\newcommand{\spinorP}{{\mathsf P}}
\newcommand{\spinorQ}{{\mathsf Q}}
\newcommand{\spinorp}{{\mathsf p}}
\newcommand{\spinorq}{{\mathsf q}}

\DeclareMathOperator{\Herm}{Herm}
\DeclareMathOperator{\trace}{tr}
\DeclareMathOperator{\diag}{diag}

\DeclareMathOperator{\SL}{\mathit{SL}}
\DeclareMathOperator{\SU}{\mathit{SU}}

\renewcommand{\phi}{\varphi}
\renewcommand{\epsilon}{\varepsilon}

\renewcommand{\theequation}{\thesection.\arabic{equation}}
\numberwithin{equation}{section}

\allowdisplaybreaks

\newtheorem{Theorem}{Theorem}
\newtheorem{Proposition}{Proposition}
\newtheorem{Corollary}{Corollary}
\newtheorem{Lemma}{Lemma}
\theoremstyle{definition}
\newtheorem{Definition}{Definition}
\theoremstyle{remark}
\newtheorem*{Remark}{Remark}

\begin{document}
\title{Hyperbolic constant mean curvature one surfaces:
Spinor representation and trinoids in hypergeometric functions}
\author{Alexander I.\,Bobenko\footnote{E--mail: {\tt bobenko@math.tu-berlin.de}}
\\
Tatyana V.\,Pavlyukevich\footnote{E--mail: {\tt
tatiana@sfb288.math.tu-berlin.de}}
\\
Boris A.\,Springborn\footnote{E--mail: {\tt springb@math.tu-berlin.de}}}
\date{Institut f\"ur Mathematik,
Technische Universit\"at Berlin, \\
Strasse des 17. Juni 136, 10623 Berlin, Germany}
\maketitle

\section{Introduction}

\begin{figure}
\centering
\mbox{\subfigure[]{\epsfig{figure=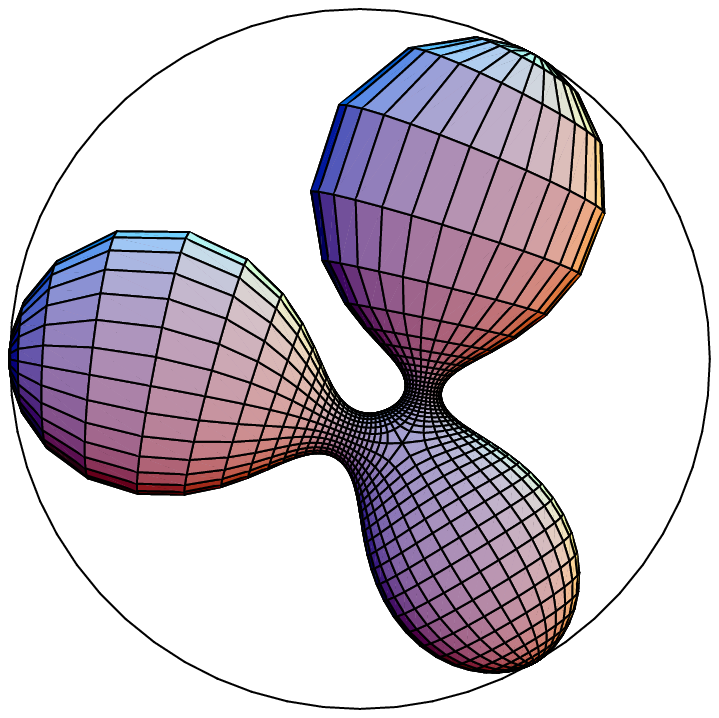,width=.40\textwidth}}\quad
    \subfigure[]{\epsfig{figure=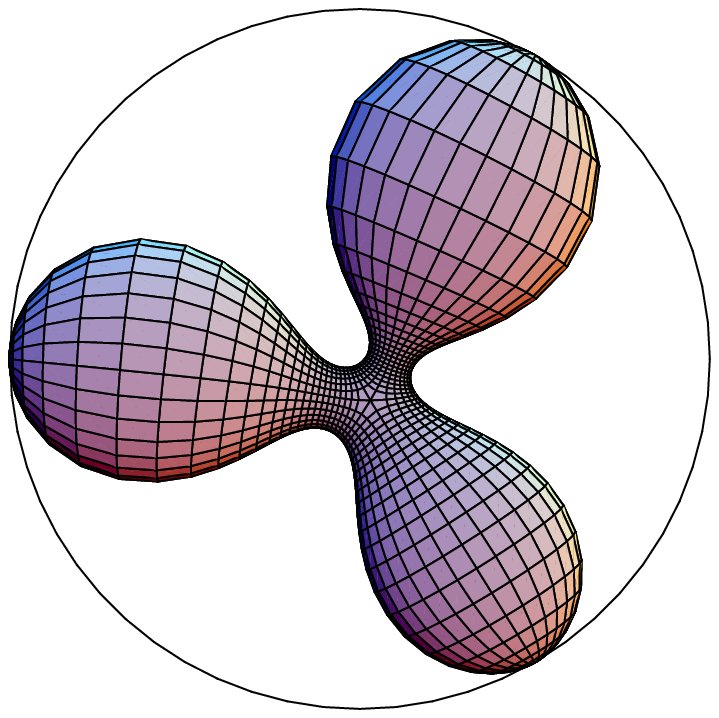,width=.40\textwidth}}}
\caption{Non-symmetric trinoids \label{f.non-symmetric}}
\end{figure}
For minimal surfaces in $\mathbb{R}^3$ there is a representation, due to
Weierstrass, in terms of holomorphic data. The Gauss-Codazzi equations for
minimal surfaces in $\mathbb{R}^3$ are equivalent to those for surfaces in
hyperbolic space with constant mean curvature 1 (CMC-1 surfaces). This lead
Bryant \cite{Bryant:1987} to derive a representation for CMC-1 surfaces in
terms of holomorphic data.

The holomorphic data used in the Weierstrass representation for minimal
surfaces consists alternatively of a function and a one-form, or of two
spinors with the same spin structure \cite{Bobenko:1994,
  Kusner/Schmitt:1996}. These functions, forms, and spinors are defined on
the same Riemann surface as the conformal minimal immersion which they
represent.  Bryant's representation for CMC-1 surfaces also involves two
spinors with the same spin structure. Other researchers prefer an equivalent
version involving a function and a one-form \cite{Umehara/Yamada:1993,
  Collin/Hauswirth/Rosenberg:2001}.  But the functions, forms, and spinors
that comprise the holomorphic data for Bryant's representation are {\em
  not}\/ defined on the same Riemann surface as the conformal immersion they
represent. As a result, a considerable amount of the great power of complex
function theory is lost. In particular, Bryant's representation does not
yield explicit formulas for CMC-1 surfaces unless their topology is very
simple.

In this paper, we present a different representation for CMC-1 surfaces in
terms of holomorphic spinors which are defined on the same Riemann surface as
the immersion. This {\em global}\/ representation is only a slight
modification of Bryant's representation, but it is much more useful if one
wants to derive explicit formulas for CMC-1 surfaces. We present a derivation
of both representations based on the method of moving frames.

We use the global representation to derive explicit formulas for CMC-1
surfaces of genus $0$ with three regular ends which are asymptotic to
catenoid cousins (CMC-1 trinoids). These surfaces were classified by Umehara
and Yamada \cite{Umehara/Yamada:2000}, but they do not present explicit
formulas.

\section{The spinor representation of surfaces in $\mathbb H^3$}

Minkowski 4-space $\mathcal L^4$ with the ca\-no\-ni\-cal
Lo\-rent\-zian met\-ric of sig\-na\-tu\-re \linebreak
$(-,+,+,+)$ can be represented as the space of
$2\times 2$ hermitian matrices. We identify $(x_0,x_1,x_2,x_3) \in
\mathcal L^4 $ with the matrix
\begin{equation*}
X=x_o I+\sum_{\alpha =1}^3 x_{\alpha}\overline{\sigma}_{\alpha}=
\begin{pmatrix}
x_0+x_3&x_1+i x_2\\
x_1-i x_2&x_0-x_3
\end{pmatrix} \in \Herm(2).  
\end{equation*}
where $\overline{\sigma}_{\alpha}$ are complex conjugate Pauli
matrices
$$
\overline{\sigma}_1=\begin{pmatrix}0&1\\ 1&0
\end{pmatrix}=\sigma_1, \overline{\sigma}_2=
\begin{pmatrix}0&i\\ -i&0\end{pmatrix}=-\sigma_2,
\overline{\sigma}_3= \begin{pmatrix}1&0\\
0&-1\end{pmatrix}=\sigma_3.
$$
In terms of the corresponding matrices the scalar product of
vectors $X$ and $Y$ is
\begin{equation*}
\langle
 X,\,Y \rangle\,=\,
 -\frac{1}{2}\;\trace(X\,\sigma_2\,Y^T\,\sigma_2).
\end{equation*}
Under this identification, hyperbolic 3-space
\begin{equation*}
\mathbb H^3 = \lbrace (x_0,x_1,x_2,x_3) \in \mathcal
L^4 ; \sum_{i=1}^3 x_i^2-x_0^2=-1, x_0>0 \rbrace
\end{equation*}
is represented as
\begin{equation*}
\begin{aligned}
\mathbb H^3
& = \lbrace X \in \Herm(2) ; \langle X,X\rangle=-1=-\det(X),
\trace(X)>0 \rbrace \\
& = \lbrace a\cdot a^*; a \in \SL(2,\mathbb C) \rbrace,
\end{aligned}
\end{equation*}
where $a^* = \overline a^T$.

Consider a smooth orientable surface in hyperbolic 3-space. The
induced metric $\Omega$ generates the complex structure of a
Riemann surface $\mathcal R$. The surface is given by an immersion
$F=(F_0,F_1,F_2,F_3):\,\mathcal R\,\to \mathbb H^3,$ and the
metric is conformal: $\Omega=e^u\,dzd\bar{z}$ where $z=x+iy$ is a local
coordinate on $\mathcal R$. The conformality of the parameterization is
equivalent to
\begin{equation*}
\langle F_z,\,F_z\rangle =
\langle F_{\overline z},\,F_{\overline z}\rangle = 0,  \quad
\langle F_z,\,F_{\overline z}\rangle =
\frac{1}{2}\,e^u.  
\end{equation*}
Here $F_z, F_{\overline z}$ are the partial derivatives with
$$
\frac{\partial}{\partial z} =
\frac{1}{2}(\frac{\partial}{\partial x}-i\,\frac{\partial}{\partial y}),\quad
\frac{\partial}{\partial\overline z} =
\frac{1}{2}(\frac{\partial}{\partial x}+i\,\frac{\partial}{\partial y}).
$$
The vectors $F, F_x, F_y$ and the unit normal $N$ define an
orthogonal moving frame on the surface
$$
\langle F,F\rangle=-1,\quad \langle N,N\rangle=1.
$$
The first  and the second fundamental forms are
\begin{align*}
\langle dF,\,dF\rangle & =  e^u\,dz\,d\overline z ,\\
-\langle dF,\,dN\rangle & =  Q\,dz^2 + H\,e^u\,dz\,d\overline z +
\overline  Q\,d\overline z^2,
\end{align*}
where
\begin{equation*}
Q = \langle F_{z\,z},\,N\rangle,\quad H\,e^u = 2\,\langle
F_{z\,\overline z},\,N\rangle. 
\end{equation*}
Here, $Q\,dz^2$ is the Hopf differential and $H$ is the mean curvature of $F$.

Conformal immersions in ${\mathbb H}^3$ can be described locally, on a domain
$D\subset {\mathbb C}$, by a smooth mapping $\phi:\,D\,\to\SL(2,\mathbb C)$
which transforms the basis $I, \overline{\sigma}_1,
\overline{\sigma}_2,\overline{\sigma}_3$ into the moving frame $F, F_x,\,
F_y,\, N$:
\begin{equation*}
\begin{aligned}
  F &= \phi\,\phi^*,\\
F_x &= e^{u/2}\,\phi\,\overline{\sigma}_1\,\phi^*,\\
F_y &= e^{u/2}\,\phi\,\overline{\sigma}_2\,\phi^*,\\
  N &= \phi\,\overline{\sigma}_3\,\phi^*.
\end{aligned}
\end{equation*}
In the complex coordinate $z=x+iy$ we have
\begin{equation*}
dF = e^{u/2}\,\phi\,
\begin{pmatrix}
0&dz\\
d\overline z&0
\end{pmatrix}
\,\phi^*.
\end{equation*}
The Gauss-Weingarten equations in terms of $\phi$ are
\begin{alignat}{2}
\phi_z &= \phi\,\widetilde U, &\qquad \widetilde U &=
\begin{pmatrix}
u_z/4&\frac{1}{2}\,(H+1)\,e^{u/2}\\
-Q\,e^{-u/2}&-u_z/4
\end{pmatrix},
\label{DE_phi_U}\\
\phi_{\overline z} &= \phi\,\widetilde V, &\qquad \widetilde V &=
\begin{pmatrix}
-u_{\overline z}/4&\overline Q\,e^{-u/2}\\
-\frac{1}{2}\,(H-1)\,e^{u/2}&u_{\overline z}/4
\end{pmatrix}.
\label{DE_phi_V}
\end{alignat}
Their compatibility condition are the Gauss-Codazzi equations
\begin{equation}                                         \label{GC}
\begin{aligned}
&u_{z\,\overline z}\,+\,\frac{1}{2}\,(H^2-1)\,e^u\,-
\,2\,Q\,\overline Q\,e^{-u} &= 0,\\
&\overline Q_z = \frac{1}{2}\,H_{\overline z}\,e^u,\\
& Q_{\overline z} = \frac{1}{2}\,H_z\,e^u.
\end{aligned}
\end{equation}

Globally, not $\phi$ but
\begin{equation}
\Phi\,
\begin{pmatrix} \sqrt{dz}&0\\
0&\sqrt{d\overline{z}}
\end{pmatrix}                                       \label{spinor}
\end{equation}
is well defined, where $\Phi =e^{u/4}\,\phi$.  This is a spinor on the
Riemann surface $\mathcal R$; it is independent of the choice of a local
coordinate $z$ on $\mathcal R$. Note that $\det \Phi =e^{u/2}$.

We arrive at the following

\begin{Theorem}                                     \label{t.spinor}
  A conformal immersion $F:{\mathcal R}\to {\mathbb H}^3$ with Gauss map $N$
  defines, uniquely up to sign, a spinor \eqref{spinor} on $\mathcal R$ such
  that locally
\begin{equation}
\begin{aligned}
  F &= e^{-u/2}\Phi\,\Phi^*,\\
  dF&= \Phi\,\begin{pmatrix}
0&dz\\
d\overline{z}&0
\end{pmatrix}
\Phi^*,\\
  N &= e^{-u/2}\Phi\,\overline{\sigma}_3\,\Phi^*.
\end{aligned}                                           \label{frame_global}
\end{equation}
Furthermore, $e^{u/2}=\det \Phi$ and
\begin{equation}                                                     \label{UV}
\begin{aligned}
\Phi^{-1}\,\Phi_z &= U, &\qquad U &=
\begin{pmatrix}
u_z/2&\frac{1}{2}\,(H+1)\,e^{u/2}\\
-Q\,e^{-u/2}&0
\end{pmatrix},\\
\Phi^{-1}\,\Phi_{\overline z} &= V, &\qquad  V &=
\begin{pmatrix}
0&\overline Q\,e^{-u/2}\\
-\frac{1}{2}\,(H-1)\,e^{u/2}&u_{\overline z}/2
\end{pmatrix}.
\end{aligned}
\end{equation}
Conversely, given a spinor \eqref{spinor} on $\mathcal R$ with $\Phi$
satisfying \eqref{UV}, where $e^{u/2}=\det \Phi$, formulas
\eqref{frame_global} describe a conformally parametrized surface in
$\mathbb H^3$ and its Gauss map $N$.
\end{Theorem}


\section{The Weierstrass representation for CMC-1 surfaces in $\mathbb H^3$}

Let $F$ be a surface in $\mathbb H^3$ with constant mean curvature $H = 1$
(CMC-1 surface). The corresponding $\Phi$ of theorem~\ref{t.spinor} satisfies
\begin{equation}
\begin{aligned}
\Phi_z &= \Phi \,U, &\qquad U &=
\begin{pmatrix}
u_z/2 &e^{u/2}\\
-Q\,e^{-u/2}&0
\end{pmatrix},\\
\Phi_{\overline z} &= \Phi \,V, &\qquad V &=
\begin{pmatrix}
0&\overline Q\,e^{-u/2}\\
0&u_{\overline z}/2
\end{pmatrix},
\label{DE_Phi_UV}
\end{aligned}
\end{equation}
Since, by the second equation, the $\overline{z}$-derivative of the first
column of $\Phi$ vanishes,
\begin{equation}\label{PQ}
  \Phi=
  \begin{pmatrix}
    \spinorP & * \\
    \spinorQ & *
  \end{pmatrix},
\end{equation}
where $\spinorP$ and $\spinorQ$ are holomorphic spinors on $\mathcal R$; see
\eqref{spinor}. Furthermore, the first equation of \eqref{DE_Phi_UV},
equation \eqref{PQ}, and $\det\Phi=e^{u/2}$ imply that the Hopf differential
is related to $\spinorP$ and $\spinorQ$ by
\begin{equation}\label{HopfPQ}
  Q = \spinorP'\spinorQ-\spinorQ'\spinorP.
\end{equation}
The hyperbolic Gauss map (see \cite{Umehara/Yamada:1996}) is
\begin{equation}\label{GaussMapPQ}
  G = -\spinorP/\spinorQ.
\end{equation}
Setting $H=1$ in \eqref{GC}, one obtains the Gauss-Codazzi equations for
CMC-1 surfaces:
\begin{equation*}
\begin{aligned}
&u_{z\,\overline z}\,-\,2\,Q\,\overline Q\,e^{-u} &= 0,\\
& Q_{\overline z} = 0.
\end{aligned}
\end{equation*}
They are invariant with respect to the transformation
\begin{equation}
\begin{aligned}
Q &\to \lambda\,Q,\\
e^{u} &\to |\lambda|^2\,e^{u}, \qquad \lambda \in \mathbb{C}\setminus\{0\}.
\end{aligned}
\label{transform}
\end{equation}
Thus, every CMC-1 surface $F$  in $\mathbb H^3$ possesses a
two-parameter family $F_\lambda$ of deformations \eqref{transform}
within the CMC-1 class.

Consider the corresponding $\Phi(z,\,\overline z,\,\lambda)$ which is a
solution of the system
\begin{alignat}{2}\label{UVoflambda}
  \Phi_z &= \Phi \,U(\lambda), &\qquad U(\lambda) &=
  \begin{pmatrix}
    u_z/2 & |\lambda| e^{u/2}\\
    -\frac{\lambda}{|\lambda|}Q\,e^{-u/2}&0
  \end{pmatrix},\\
  \Phi_{\overline z} &= \Phi \,V(\lambda), &\qquad V(\lambda) &=
  \begin{pmatrix}
    0&\frac{\overline{\lambda}}{|\lambda|}\overline Q\,e^{-u/2}\\
    0&u_{\overline z}/2
  \end{pmatrix}.
\end{alignat}
Now let $\lambda\rightarrow 0$ while $\frac{\lambda}{|\lambda|}=1$. The
corresponding equations have solutions of the form
\begin{equation}
\Phi_0=
\begin{pmatrix}
 \spinorp & \overline \spinorq\\
-\spinorq & \overline \spinorp
\end{pmatrix}              \label{pq}
\end{equation}
where $\spinorp$ and $\spinorq$ are holomorphic spinors on
the universal covering $\widetilde {\mathcal R}$ of $\mathcal R$, and
\begin{equation}\label{metric&hopf}
  \begin{split}
    e^{u/2} &= |\spinorp|^2+|\spinorq|^2,\\
          Q &= -\spinorp'\spinorq+\spinorp\spinorq'.
  \end{split}
\end{equation}

\begin{Remark}
  Note that $\spinorP$ and $\spinorQ$ are well defined holomorphic spinors on
  the Riemann surface $\mathcal R$, but the spinors $\spinorp$ and $\spinorq$
  are only well defined on the universal cover $\widetilde{\mathcal R}$ of
  $\mathcal R$.
\end{Remark}

Let $\Phi_1=\Phi|_{\lambda=1}$ and denote by $\Psi$ the quotient
\begin{equation}
\Phi_1\,=\,\Psi\,\Phi_0\,.                                  \label{Psi}
\end{equation}
\begin{Theorem}
  The mapping $\Psi:\widetilde {\mathcal R}\to \SL(2,\mathbb C)$ defined by
  \eqref{Psi} is holomorphic and satisfies
  \begin{align}
    \label{WR_pq}
    \Psi_z &= \Psi
    \begin{pmatrix}
      \spinorp\spinorq & \spinorp^2\\
      -\spinorq^2      & -\spinorp\spinorq
    \end{pmatrix},\\
    \label{WR_PQ}
    \Psi_z &=
    \begin{pmatrix}
      \spinorP\spinorQ & \spinorP^2\\
      -\spinorQ^2      & -\spinorP\spinorQ
    \end{pmatrix}
    \Psi,
  \end{align}
  where $\spinorp$, $\spinorq$ are the holomorphic spinors on
  $\widetilde {\mathcal R}$ defined by \eqref{pq}, and $\spinorP$,
  $\spinorQ$ are the holomorphic spinors on $\mathcal R$ defined by
  \eqref{PQ}.

  The immersion $F:\mathcal R\to \mathbb H^3$ is recovered by
  \begin{equation}\label{immersionformula}
    F=\Psi\Psi^*.
  \end{equation}
\end{Theorem}

\begin{proof}  
  Since
  \begin{equation*}
    \Psi_{\overline z}=(\Phi_1\,{\Phi_0}^{-1})_{\overline z}
    =\Phi_1({\Phi_1}^{-1}\Phi_{1\,\overline{z}}
            - {\Phi_0}^{-1}\Phi_{0\,\overline{z}}){\Phi_0}^{-1},
  \end{equation*}
  equations \eqref{UVoflambda} imply $\Psi_{\overline z}=0$. Hence $\Psi$ is
  holomorphic.

  Similarly one finds that,
  \begin{equation*}
    \Psi_z = e^{u/2}\Phi_1
    \begin{pmatrix}
      0 & 1\\
      0 & 0
    \end{pmatrix}
    \Phi_0^{-1},
  \end{equation*}
  and hence,
  \begin{equation*}
    \Psi^{-1}\,\Psi_z\,=\,e^{u/2}\,\Phi_0\,
    \begin{pmatrix}
      0&1\\
      0&0
    \end{pmatrix}\,\Phi_0^{-1}
  \end{equation*}
  and
  \begin{equation*}
    \Psi_z\,\Psi^{-1}\,=\,e^{u/2}\,\Phi_1\,
    \begin{pmatrix}
      0&1\\
      0&0
    \end{pmatrix}\,\Phi_1^{-1}.
  \end{equation*}
  Now, \eqref{pq} and \eqref{PQ} imply \eqref{WR_pq}, \eqref{WR_PQ}.
  Finally, equation \eqref{frame_global} and $\Phi_0{\Phi_0}^* = e^{u/2}I$
  imply the immersion formula \eqref{immersionformula}.
\end{proof} 

By equation \eqref{WR_pq} and the immersion formula \eqref{immersionformula},
the spinors $\spinorp$ and $\spinorq$ determine the surface
$F$ up to a hyperbolic isometry. The metric and Hopf differential are related
to $\spinorp$ and $\spinorq$ by \eqref{metric&hopf}. This representation of
CMC-1 surfaces, which is due to Bryant \cite{Bryant:1987}, is therefore an
{\em intrinsic} and {\em metric} description. It is also essentially {\em
  local}, since the spinors $\spinorp$ and $\spinorq$ are
not well defined on the Riemann surface $\mathcal{R}$, but only on its
universal cover. This is a serious disadvantage if one wants to construct
CMC-1 surfaces with non-trivial topology. In particular, this prohibits in
all but the simplest cases the integration of equation \eqref{WR_pq} in
closed form.

Formula \eqref{WR_PQ}, on the other hand, is a {\em global}\/ representation
of a CMC-1 surface by holomorphic spinors $\spinorP$ and
$\spinorQ$ on $\mathcal{R}$. Unfortunately, these spinors do in
general not determine the surface up to isometry. While the Hopf
differential and the hyperbolic Gauss map are determined by \eqref{HopfPQ}
and \eqref{GaussMapPQ}, the metric depends non-trivially on the particular
solution of \eqref{WR_PQ}. But there is also the global condition that the
immersion $F$ obtained from \eqref{immersionformula} is well defined on
$\mathcal{R}$. This, together with $\spinorP$ and
$\spinorQ$, may determine the surface uniquely if $\mathcal{R}$ is
not simply connected. The condition that $F$ is well defined on $\mathcal{R}$
implies the following corollary.


\begin{Corollary}
  Let $F:\mathcal R\to\mathbb H^3$ be a CMC-1 surface in $\mathbb H^3$ and
  $\Phi$ its spinor frame \eqref{PQ}, defining holomorphic spinors
  $\spinorP$ and\/ $\spinorQ$ on $\mathcal R$. Then equation
  \eqref{WR_PQ}
  has a solution $\Psi: \widetilde {\mathcal R}\to \SL(2,\mathbb C)$ with
  unitary monodromy.
\end{Corollary}

Conversely, one obtains the following representation theorem.

\begin{Theorem}                                 \label{t.unitary}
  Let $\spinorP$ and\/ $\spinorQ$ be two holomorphic
  spinors with the same spin structure on a Riemann surface $\mathcal R$ and
  suppose $\Psi: \widetilde {\mathcal R}\to \SL(2,\mathbb C)$ a solution of
  equation \eqref{WR_PQ}
  with unitary monodromy. Then equation \eqref{immersionformula}
  defines a CMC-1 immersion $F:{\mathcal R}\to\mathbb H^3$.
\end{Theorem}

Rossman, Umehara, Yamada, and others describe CMC-1 surfaces in terms of the
`secondary Gauss map' $g = -\spinorp/\spinorq$ and the one-form $\omega =
-\spinorq^2\,dz$. Thus, instead of equation \eqref{WR_pq}, they write
\begin{equation*}
  d\Psi = \Psi
  \begin{pmatrix}
    g & -g^2 \\
    1 & -g
  \end{pmatrix}\omega.
\end{equation*}
The secondary Gauss map $g$ and one-form $\omega$ are not defined on the same
Riemann surface $\mathcal R$ on which the conformal immersion $F$ is defined.
The hyperbolic Gauss map $G=-\spinorP/\spinorQ$ and the holomorphic one-form
$\Omega=-\spinorQ^2\,dz$, one the other hand, are defined on the Riemann
surface $\mathcal R$. In terms of these, equation \eqref{WR_PQ} reads
\begin{equation*}
  d\Psi=
  \begin{pmatrix}
    G & -G^2 \\
    1 & -G
  \end{pmatrix}\Psi\Omega.
\end{equation*}
Even though Rossman, Umehara and Yamada are aware of this equation
\cite{Rossman/Umehara/Yamada:1997}, they do not consider $G$ and $\Omega$ as
the Weierstrass data for the CMC-1 immersion $F$ but for a dual immersion.

\section{Catenoid cousin, catenoidal ends and n-noids}\label{s.catenoid}

\begin{figure}[ht]
\centering
\mbox{\subfigure[$\lambda<\tfrac{1}{2}$]{\epsfig{figure=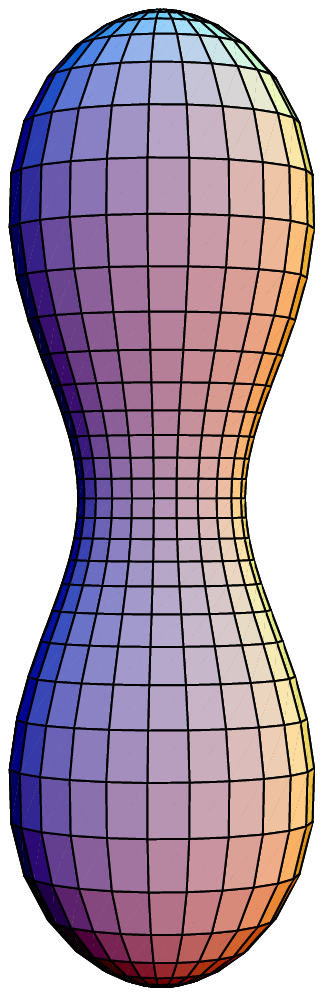,width=.40\textwidth}}\quad
    \subfigure[$\lambda>\tfrac{1}{2}$]{\epsfig{figure=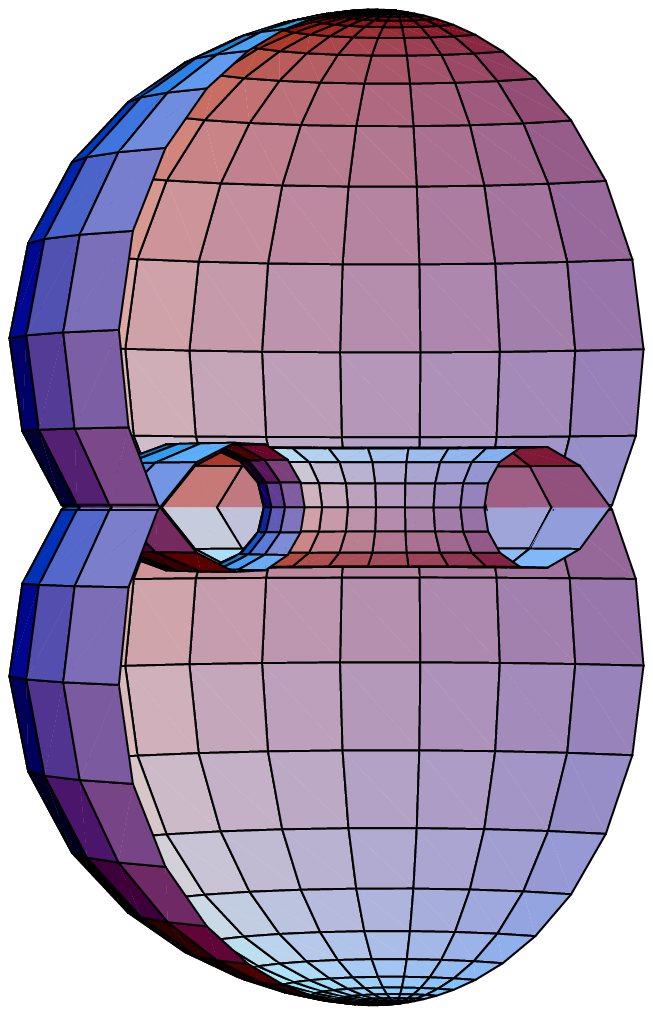,width=.40\textwidth}}}
\caption{CMC-1 twonoids in the Poincar\'{e} model of $\mathbb H^3$.
\label{twonoids}}
\end{figure}

Let us start our investigation of special CMC-1 surfaces in $\mathbb H^3$
with a simple example of the catenoid cousins which we also call
\emph{twonoids}. Since the Gauss equations of CMC-1 surfaces in $\mathbb H^3$
and of minimal surfaces in $\mathbb R^3$ coincide these surfaces are locally
isometric. The catenoid cousins are surfaces isometric to the catenoids. They
were investigated by Bryant \cite{Bryant:1987}.

These surfaces are of genus zero with two regular ends. In our global
spinorial description, twonoids are immersions
\begin{equation*}
F=\Psi\Psi^*:\mathbb C\setminus\{0\}\to\mathbb{H}^3,
\end{equation*}
where $\Psi$ satisfies the differential equation \eqref{WR_PQ}
with the Weierstrass data
\begin{equation*}
\spinorP=\frac{p_0}{z}+p_{\infty}, \qquad \spinorQ=\frac{q_0}{z}+q_{\infty}.
\end{equation*}
(By applying a suitable hyperbolic isometry and a coordinate transformation
$z\to az$ to $F$ one can reduce this to the simpler case $p_0=q_{\infty}=0$,
$p_{\infty}=q_0$.)  This equation can be solved explicitly in elementary
functions. A particular solution with determinant 1 is
\begin{equation*}
  \Psi_0 = cB
  \begin{pmatrix}
    z^{1/2} & 0 \\
    0       & z^{-1/2}
  \end{pmatrix} C
  \begin{pmatrix}
    z^{\lambda} & 0\\
    0         & z^{-\lambda}
  \end{pmatrix},
\end{equation*}
where
\begin{align*}
  B&=
  \begin{pmatrix}
     \left(\frac{p_0}{z}+p_{\infty}\right)
              & \frac{p_0}{p_0 q_\infty-p_\infty q_0}\\
    -\left(\frac{q_0}{z}+q_{\infty}\right)
              & -\frac{q_0}{p_0 q_\infty-p_\infty q_0}
  \end{pmatrix},\\
  C&=
  \begin{pmatrix}
    \frac{2\lambda-1}{2(p_0 q_\infty-p_\infty q_0)}
               & -\frac{2\lambda+1}{2(p_0 q_\infty-p_\infty q_0)}\\
    1          & 1
  \end{pmatrix},
\end{align*}
\begin{equation*}
  \lambda = \frac{1}{2}\sqrt{1+4(p_0 q_\infty-p_\infty q_0)},\quad
  c = \sqrt{\frac{p_0 q_\infty-p_\infty q_0}{2\lambda}}.
\end{equation*}
The general solution with determinant 1 is $\Psi=\Psi_0 A$, with $A\in
\SL(2,\mathbb C)$. Since multiplying $A$ on the right with a unitary matrix
does not change the immersion $F$, we may assume $A$ to be hermitian. When
continued along a path going around the puncture $z=0$ in the
counterclockwise direction, $\Psi$ is transformed into $\Psi\mathcal{M}_0$,
where the monodromy matrix is
\begin{equation*}
  \mathcal M_0 = -A^{-1}
  \begin{pmatrix}
    e^{2\pi i\lambda} & 0 \\
    0                 & e^{-2\pi i\lambda}
  \end{pmatrix} A.
\end{equation*}
For $\mathcal{M}_0$ to be unitary, $\lambda$ must be real. If $\lambda$ is
not half-integer, then $A$ must be diagonal. In fact, it suffices to consider
$A=I$, since different $A$ yield the same surface up to a hyperbolic isometry
and a coordinate change $z\to az$.
If $\lambda$ is half-integer, then $A$ is arbitrary. In this case, one
obtains also surfaces which are not surfaces of revolution, and which are not
locally isometric to a catenoid.

For the surfaces of revolution, the profile curve is embedded if
$\lambda<\frac{1}{2}$, and it has a single self-intersection if
$\lambda>\frac{1}{2}$, see Fig.~\ref{twonoids}.


There are no compact CMC-1 surfaces in $\mathbb H^3$. Bryant has shown that
the Riemann surface of a complete conformal immersion $F:\mathcal R\to
\mathbb H^3$ of finite total curvature can be compactified: $\mathcal
R=\hat{\mathcal R}\setminus \left\{a_1,a_2,\dotsc,a_N\right\}$, where
$\hat{\mathcal R}$\/ is a compact Riemann surface \cite{Bryant:1987}.
Moreover, Collin, Hauswirth, and Rosenberg have shown that a properly
embedded annular end is of finite total curvature and regular
\cite{Collin/Hauswirth/Rosenberg:2001}. The punctures $a_1,a_2,\dotsc,a_N$
correspond to the ends of the immersion.  For their classification one uses
the hyperbolic Gauss map $G=-\spinorP/\spinorQ$.
The end corresponding to a point $a_i\in \hat{\mathcal R}$ is called
\emph{regular} if $G$ can be meromorphically extended to $a_i$, and
\emph{irregular} if it is an essential singularity of $G$. Motivated by the
behavior of the Weierstrass data at the punctures of twonoids, it is natural
to give the following analytic definition of the catenoidal ends.

\begin{Definition}
  The end corresponding to a puncture $a_i$, is called \emph{catenoidal}, if
  the spinors $\spinorP$, $\spinorQ$ have only simple poles
  at $a_i$.
\end{Definition}

I.~e., it is required that, for a local coordinate $z$
centered in $a_i$, the Weierstrass data $\spinorP$, $\spinorQ$ satisfy
\begin{equation*}
  \spinorP=\frac{p_0}{z}+O(1) \quad\text{and}\quad
  \spinorQ=\frac{q_0}{z}+O(1)
  \quad \text{for}\quad z\to 0.
\end{equation*}
Obviously, catenoidal ends are regular.

We call a compact CMC-1 surface of genus zero with $n$ catenoidal ends an
\emph{$n$-noid}. Normalizing one end to $z=\infty$, all $n$-noids can be
conformally parametrized as
\begin{equation*}
  F:\mathbb{C}\setminus\{a_1,a_2,\dotsc,a_{n-1}\}\to\mathbb H^3
\end{equation*}
with the Weierstrass data
\begin{equation}                              \label{Wei-n-noids}
\spinorP=\sum_{i=1}^{N-1}\frac{p_i}{z-z_i} + p_{\infty}, \quad
\spinorQ=\sum_{i=1}^{N-1}\frac{q_i}{z-z_i} + q_{\infty}.
\end{equation}
%
%
At a catenoidal end, the system \eqref{WR_PQ} is locally gauge equivalent to
a Fuchsian system. Indeed, let $z=0$ be a puncture and suppose
$\spinorP$ and $\spinorQ$ satisfy
\begin{equation}                        \label{PQ-local}
\spinorP=\frac{a_{-1}}{z}+a_0+o(1), \quad \spinorQ=\frac{b_{-1}}{z}+b_0+o(1)
\quad\text{for}\quad z\to 0.
\end{equation}
The following lemma is obtained by direct calculation.
\begin{Lemma}
  If $\Psi$ satisfies equation \eqref{WR_PQ} with $\spinorP$, $\spinorQ$ as
  in \eqref{PQ-local}, then the gauge equivalent $\widetilde{\Psi}$ defined
  by
  \begin{equation*}
    \Psi=
    \begin{pmatrix}
      a_{-1}&0\\
      -b_{-1}&\frac{1}{a_{-1}}
    \end{pmatrix}\,
    \begin{pmatrix}
      \frac{1}{\sqrt z}&0\\
      0&\sqrt z
    \end{pmatrix}
    \widetilde \Psi
  \end{equation*}
  satisfies an equation $\widetilde{\Psi}_z=\widetilde{A}\widetilde{\Psi}$,
  with
  \begin{equation*}
    \widetilde A = \frac{1}{z}
    \begin{pmatrix}
      \frac{1}{2} + r & 1 \\
      -r^2            & -\frac{1}{2} - r
    \end{pmatrix}+O(1)\quad\text{for}\quad z\to\infty,
  \end{equation*}
  where $r = a_{-1} b_0 - a_0 b_{-1}$.
\end{Lemma}

\begin{Corollary}\label{localMonodromy}
  Under the conditions of the lemma, the local monodromy of \eqref{WR_PQ}
  around $z=0$ is
  \begin{equation*}
    M =
    \begin{pmatrix}
      e^{2\pi i\alpha} & 0 \\
      0                   & e^{-2\pi i\alpha}
    \end{pmatrix}
    \quad\text{with}\quad
    \alpha = \frac{1}{2} + \sqrt{\frac{1}{4} + r}.
  \end{equation*}
\end{Corollary}



\section{Trinoids. Reduction to a Fuchsian system}

The rest of the paper is devoted to explicit description of the
\emph{trinoids}, which are CMC-1 immersions of genus zero with three
catenoidal ends. Without loss of generality the punctures can be normalized
to $0,1,\infty$. By equation \eqref{Wei-n-noids}, the trinoids are thus
conformal immersions $F:\mathbb{C}\setminus\{0,1\}\to\mathbb{H}^3$
with Weierstrass data
\begin{equation}                   \label{W_data}
\spinorP = \frac{p_0}{z} + \frac{p_1}{z-1} + p_{\infty},\quad
\spinorQ = \frac{q_0}{z} + \frac{q_1}{z-1} + q_{\infty}.
\end{equation}
The asymptotics at $z=0$, $z=1$ and $z=\infty$ are as follows.
\begin{align*}
  &z\to 0 :
  &&\spinorP = \frac{p_0}{z} + (p_\infty-p_1) + o(1),
  &&\spinorQ = \frac{q_0}{z} + (q_\infty-q_1) + o(1), \\
  &z\to 1 :
  &&\spinorP = \frac{p_1}{z-1}+(p_0+p_\infty)+o(1),
  &&\spinorQ = \frac{q_1}{z-1}+(q_0+q_\infty)+o(1), \\
  &z\to\infty :
  &&\spinorP = p_\infty + \frac{p_0 + p_1}{z} + o(1),
  &&\spinorQ = q_\infty + \frac{q_0 + q_1}{z} + o(1).
\end{align*}
By corollary \ref{localMonodromy}, the the local monodromy around
$j=0,1,\infty$ is
\begin{equation}\label{loc_mon}
  M_j=
  \begin{pmatrix}
    e^{2\pi i \alpha_z}&0\\
    0&e^{-2\pi i \alpha_z}
  \end{pmatrix},
  \qquad\alpha_j = \frac{1}{2}+\sqrt{\frac{1}{4}+c_j},
\end{equation}
where
\begin{equation}   \label{c}
  \begin{aligned}
    c_0      &= \langle p,q\rangle_{10}+\langle p,q\rangle_{0\infty},\\
    c_1      &= \langle p,q\rangle_{10}+\langle p,q\rangle_{1\infty},\\
    c_\infty &= \langle p,q\rangle_{0\infty}+\langle p,q\rangle_{1\infty},
  \end{aligned}
\end{equation}
and
\begin{equation*}
\langle p,q\rangle_{ij} = p_i q_j - p_j q_i,\quad
i\neq j,\quad i,j = 0,1,\infty.
\end{equation*}

In our integration of trinoids we proceed as follows. First, we show that the
corresponding system \eqref{WR_PQ} is globally gauge equivalent to a Fuchsian
system with three singularities. The latter can be solved explicitly in terms
of hypergeometric functions. This provides explicit formulas for the
monodromy matrices of the original system. By theorem~\ref{t.unitary},
trinoids are obtained if the monodromy matrices are unitary.

\begin{Proposition}    \label{th_gauge}
  If $\Psi$ satisfies equation \eqref{WR_PQ} with $\spinorP$, $\spinorQ$ as
  in equation \eqref{W_data}, then $\Phi$ defined by $\Psi=D\Phi$,
  \begin{equation}     \label{gauge}
    D=
    \begin{pmatrix}
      P&\alpha_1\,z+\beta_1\\
      -Q&\alpha_2\,z+\beta_2
    \end{pmatrix}\,
    \begin{pmatrix}
      \sqrt{z-1}&0\\
      \frac{k}{z\,\sqrt{z-1}}&\frac{1}{\sqrt{z-1}}
    \end{pmatrix}\,
    \begin{pmatrix}
      \frac{2\,\alpha}{\mu}&0\\
      1&1
    \end{pmatrix},
  \end{equation}
  satisfies the Fuchsian system
  \begin{equation}    \label{sys_Fuchs}
    \Phi_z\,=\,\left(\frac{A_0}{z}+\frac{A_1}{z-1}\right)\,\Phi,\,
  \end{equation}
  with
  \begin{equation*}
    A_0=
    \begin{pmatrix}
      \alpha&0\\
      0&-\alpha
    \end{pmatrix}\,
    A_1=
    \begin{pmatrix}
      \beta&\gamma\\
      \delta&-\beta
    \end{pmatrix}.
  \end{equation*}
  Here, the coefficients are as follows:
  \begin{equation}      \label{gauge-constants}
    \begin{aligned}
      \alpha\,&=\,\frac{1}{2}\,
      \left(1-\sqrt{1+4\,\zinf+4\,\oz}\,
      \right),\\
      \beta\,&=\,\frac{1}{2}\,\frac{\oz\,
        (1-2\,\alpha)-\zinf}
      {\zinf+\oz}\,,\\
      \gamma\,&=\zinf\,
      \left(\frac{\oinf}{\Delta}+\frac{1}{\alpha}\right)\,,\\
      \delta\,&=\,\frac{\Delta}{\zinf}\,\frac{\Delta+\oinf\,\alpha}
      {\Delta-\zinf\,\oz+(\Delta+\oinf)\,\alpha}\,,\\
      \mu\,&=\,2\,\zinf
      \left(1-k\,\frac{\oinf}{\Delta}\right)\,,\\
      k\,&=\,\Delta\,\frac{\zinf\,\oz-\Delta\,\alpha}{\Delta^2+\oz\,\zinf\,\oinf}\,,\\
      \alpha_1\,&=\,-\frac{p_\infty\,\oz}{\Delta}\,,\quad
      \alpha_2\,=\,\frac{q_\infty\,\oz}{\Delta}\,,\\
      \beta_1\,&=\,\frac{p_0\,\oinf}{\Delta}\,,\qquad
      \beta_2\,=\,-\frac{q_0\,\oinf}{\Delta}\,,\\
      \Delta\,&=\,\oz\,\zinf+\oz\,\oinf+\zinf\,\oinf\,.
    \end{aligned}
  \end{equation}
\end{Proposition}
\begin{proof} 
  We will construct the gauge transformation as a composition of three more
  elementary transformations $D=BCM$. Only the $B$ part is non-trivial.
  Construct a matrix $B=\begin{pmatrix}P& S\\-Q& T\end{pmatrix}$ with $\det
  B=1$ which transforms $A$ to its Jordan form:
  \begin{equation*}
    A\,=\,B\,\begin{pmatrix} 0&1\\0&0
    \end{pmatrix}\,B^{-1}.
  \end{equation*}
  The determinant condition can be satisfied by choosing
  \begin{equation*}
    S=\alpha_1\,z+\beta_1,\qquad T=\alpha_2\,z+\beta_2.
  \end{equation*}
  Then the condition $\det B=1$ implies the system of linear
  equations
  \begin{equation*}
    \mathcal A\,
    \begin{pmatrix}
      \alpha_1\\
      \alpha_2\\
      \beta_1\\
      \beta_2
    \end{pmatrix}\,=\,\begin{pmatrix}0\\0\\0\\1\end{pmatrix},\qquad
    \mathcal A\,=\,
    \begin{pmatrix}
      q_\infty&p_\infty&0&0\\
      0&0&q_0&p_0\\
      q_1&p_1&q_1&p_1\\
      q_0+q_1&p_0+p_1&q_\infty&p_\infty
    \end{pmatrix}\,.
  \end{equation*}
  Note that $\det \mathcal A=\Delta$. Formulas
  \eqref{gauge-constants} for
  $\alpha_1,\,\alpha_2,\,\beta_1$ and $\beta_2$ give the solution
  of the system.

  After this first gauge transformation $\Psi=B\/\widetilde \Psi$ we obtain
  the equation $\widetilde \Psi_z=\widetilde A\,\widetilde \Psi$ with
  \begin{equation*}
    \widetilde A=
    \begin{pmatrix}
      \frac{\oz\zinf}{\Delta}\frac{1}{z}+
      \frac{\oz\oinf}{\Delta}\frac{1}{z-1}&
      1+\frac{\oz\zinf\oinf}{\Delta^2}\\
      \frac{\zinf}{z^2}+\frac{\oinf}{(z-1)^2}+\frac{\oz}{z^2\,(z-1)^2}&
      -\frac{\oz\zinf}{\Delta}\frac{1}{z}-
      \frac{\oz\oinf}{\Delta}\frac{1}{z-1}
    \end{pmatrix}.
  \end{equation*}
  The next transformation
  $\widetilde \Psi=C\,\widehat \Psi$,
  $C=\begin{pmatrix}
    \sqrt{z-1}&0\\
    \frac{k}{z\,\sqrt{z-1}}&\frac{1}{\sqrt{z-1}}
  \end{pmatrix}
  $ almost brings the equation to Fuchsian form:
  $\widehat \Psi_z=
  \widehat A\,\widehat \Psi$, where
  \begin{equation*}
    \widehat A
    =\frac{1}{\Delta^2}
    \begin{pmatrix}
      \frac{\hat a_{11}^0}{z}+\frac{\hat a_{11}^1}{z-1}&\frac{\hat a_{12}^1}{z-1}\\
      \frac{\hat a_{21}^0}{z}+\frac{\hat a_{21}^1}{z-1}+\frac{\hat a_{21}^2}{z^2}&
      \frac{\hat a_{22}^0}{z}+\frac{\hat a_{22}^1}{z-1}
    \end{pmatrix},
  \end{equation*}
  where
  \begin{align*}
    \hat a_{11}^0&=-\hat a_{22}^0\,=\,\Delta\,\oz\,\zinf-k\,\hat a\,,\\
    \hat a_{11}^1&=-\hat a_{22}^1\,=
    \,\Delta\,\oz\,\oinf-\frac{\Delta^2}{2}+k\,\hat a\,,\\
    \hat a_{12}^1&=\hat a\,,\\
    \hat a_{21}^0&=\Delta^2\,(\zinf-\oz)-k\,(\Delta^2-2\,\Delta\,\oz\,\oinf)
    +k^2\,\hat a\,,\\
    \hat a_{21}^1&=\Delta^2\,(\oinf+\oz)+k\,(\Delta^2-2\,\Delta\,\oz\,\oinf)
    -k^2\,\hat a\,,\\
    \hat a_{21}^2&=-\Delta^2\,(\zinf+\oz)+k\,(\Delta^2-2\,\Delta\,\oz\,\zinf)
    +k^2\,\hat a\,,\\
    \hat a&=\Delta^2+\oz\,\zinf\,\oinf\,.
  \end{align*}
  Choosing
  \begin{equation*}
    k=\,\Delta\,\frac{\oz\,\zinf-
      \frac{1}{2}\,(1-\sqrt{1+4\,\oz+4\,\zinf})\,\Delta}
    {\Delta^2+\oz\,\zinf\,\oinf},
  \end{equation*}
  we bring $\widehat A$ to the Fuchsian form ($\hat a_{21}^2=0$):
\begin{equation*}
\widehat A\,=\,\frac{\widehat A_0}{z}+\frac{\widehat A_1}{z-1}\,,\qquad
\widehat A_0\,=\,
\begin{pmatrix}
  \alpha&0\\
  \mu&-\alpha
\end{pmatrix},\quad
\widehat A_1\,=\,
\begin{pmatrix}
  \hat \beta&\hat \gamma\\
  \hat \delta&-\hat \beta
\end{pmatrix},
\end{equation*}
with $\alpha$ and $\mu$ given by \eqref{gauge-constants} and
\begin{align*}
  \hat \beta\,&=\,
  -\frac{\zinf\,\oinf}{\Delta}+\frac{1}{2}-\alpha\,,\\
  \hat \gamma\,&=\,
  \frac{\Delta^2+\oz\,\zinf\,\oinf}{\Delta^2}\,,\\
  \hat \delta\,&=\, \frac{2\,k\,\zinf\,\oinf}{\Delta}-\zinf+\oinf.
\end{align*}
Finally, the transformation $\widehat \Psi=M\Phi$ with $ M=\begin{pmatrix}
  \frac{2\,\alpha}{\mu}&0\\
  1&1
\end{pmatrix}
$ implies \eqref{sys_Fuchs} with $\beta$, $\gamma$, $\delta$ as in
\eqref{gauge-constants}.
\end{proof} 

\section{Trinoids. Solution of the Fuchsian system}
A Fuchsian system of two first-order differential equations
with three singularities can be solved explicitly in terms of
hypergeometric functions. Let us diagonalize the singularities of
$A$:
\begin{equation}
A_0=L_0 \Lambda_0 L_0^{-1},\ A_1=L_1 \Lambda_1 L_1^{-1},\
-A_0-A_1=L_\infty \Lambda_\infty L_\infty^{-1},
\end{equation}
where
\begin{equation}
\begin{aligned}
\Lambda_0=\alpha \sigma_3,\quad
\Lambda_1=\tau \sigma_3,\quad
\Lambda_\infty=\rho \sigma_3,\\
\tau\,=\,\sqrt{\beta^2+\gamma\,\delta},\quad
\rho\,=\,\sqrt{(\alpha+\beta)^2+\gamma\,\delta}.
\end{aligned}
\end{equation}
\begin{Remark}
  For simplicity we consider in this paper only the \emph{generic} case when
  the differences of the eigenvalues of the singularities of the Fuchsian
  system are non-integer, i.~e.\
  \begin{equation}                       \label{generic}
    2\alpha, 2\tau, 2\rho \notin {\mathbb Z}.
  \end{equation}
  The case of half-integer $\alpha, \tau$ or $\rho$ can be treated similarly,
  although the computations are involved because many degenerated cases have
  to be considered considered.
\end{Remark}

Denote by $\Phi^{(0)}, \Phi^{(1)}$ and $\Phi^{(\infty)}$ the
canonical solutions of \eqref{sys_Fuchs} determined by their
asymptotics at the singularities
\begin{equation}                                            \label{can_sol}
\begin{aligned}
\Phi^{(0)}&=(L_0+o(z))z^{\Lambda_0},\qquad &z\to 0,\\
\Phi^{(1)}&=(L_1+o(z-1))(z-1)^{\Lambda_1},\qquad &z\to 1,\\
\Phi^{(\infty)}&=(L_\infty+o(1/z))z^{-\Lambda_\infty},\qquad &z\to \infty.
\end{aligned}
\end{equation}
\begin{Theorem}                                 \label{th_Fuchs}
The canonical solutions of the Fuchsian system \eqref{sys_Fuchs}
are given by
\begin{align}
&\Phi^{(0)}(z)\,=\,
\notag\\
&\begin{pmatrix}
\begin{aligned}
&-\tfrac{2\alpha+1}{\delta}z^{\alpha}\,(z-1)^{\tau}\,\cdot\\
&\F(a,b;c;z)
\end{aligned} &
\begin{aligned}
&z^{1-\alpha}\,(z-1)^{\tau}\,\cdot\\
&\F(a-c+1,b-c+1;2-c;z)
\end{aligned}\\
&\\
\begin{aligned}
&z^{1+\alpha}\,(z-1)^{\tau}\,\cdot \\
&\F(a+1,b+1;c+2;z)
\end{aligned}&
\begin{aligned}
&\tfrac{2\alpha-1}{\gamma}z^{-\alpha}\,(z-1)^{\tau}\,\cdot\\
&\F(a-c,b-c;-c;z)
\end{aligned}
\end{pmatrix}\,,
\label{Phi_0}\\
&\Phi^{(1)}(z)\,=\,
\notag\\
&\begin{pmatrix}
\begin{aligned}
&\tfrac{\beta+\tau}{\delta}
z^{\alpha}\,(z-1)^{\tau}\,\cdot\\
&\F(a,b;a+b-c+1;1-z)\\
\end{aligned} &
\begin{aligned}
&z^{\alpha}\,(z-1)^{-\tau}\,\cdot\\
&\F(c-a,c-b;c-a-b+1;1-z)
\end{aligned}\\
&\\
\begin{aligned}
&z^{-\alpha}\,(z-1)^{\tau}\,\cdot\\
&\F(a-c,b-c;a+b-c+1;1-z)
\end{aligned}&
\begin{aligned}
&-\tfrac{\beta+\tau}{\gamma}
z^{-\alpha}\,(z-1)^{-\tau}\,\cdot\\
&\F(-a,-b;c-a-b+1;1-z)
\end{aligned}
\end{pmatrix}\,,
\label{Phi_1}\\
&\Phi^{(\infty)}(z)=
\notag\\
&\begin{pmatrix}
\begin{aligned}
&\tfrac{\gamma(c-a)}{a(\beta+\tau)}
z^{-\tau-\rho}(z-1)^{\tau}\cdot\\
&\F(a,a-c+1;a-b+1;\frac{1}{z})\\
\end{aligned} &
\begin{aligned}
&z^{-\tau+\rho}(z-1)^{\tau}\cdot\\
&\F(b,b-c+1;b-a+1;\frac{1}{z})
\end{aligned}\\
&\\
\begin{aligned}
&z^{-\tau-\rho}(z-1)^{\tau}\cdot\\
&\F(a+1,a-c;a-b+1;\frac{1}{z})
\end{aligned}&
\begin{aligned}
&\tfrac{b(\beta+\tau)}{\gamma(c-b)}
z^{-\tau+\rho}(z-1)^{\tau}\cdot\\
&\F(b+1,b-c;b-a+1;\frac{1}{z})
\end{aligned}
\end{pmatrix},
\label{Phi_inf}
\end{align}
where $\F(a,b,c;z)$ is the hypergeometric function and
\begin{gather*}
a=\alpha+\tau+\rho,\quad b=\alpha+\tau-\rho,\quad c=2\alpha.
\end{gather*}
\end{Theorem}

The proof is given in Appendix~\ref{proofs}. It is a direct but long
computation.  The canonical solutions \eqref{Phi_0}--\eqref{Phi_inf} have
branch points at $z=0$, $z=1$ and $z=\infty$. We choose the branch cuts from
$1$ to $\infty$ along the positive real axis and from $0$ to $\infty$ along
the negative real axis.

Let us compute the monodromy group of system \eqref{sys_Fuchs}.
Fix a base point $a \in \widehat{\mathbb C}\setminus
\{0,1,\infty\}$ and a matrix $R_0\in \SL(2,\mathbb C)$. Let
$\Phi(z)$ be a solution of \eqref{sys_Fuchs} with $\Phi(a)=R_0$.
Its analytic continuation $\Phi_\gamma(z)$ along a loop
$\gamma\,\in\,\pi_1(\widehat{\mathbb C}\setminus \{0,1,\infty\})$ determines
the monodromy matrix $M_\gamma\in \SL(2,\mathbb C)$ through
\begin{equation*}
\Phi_\gamma(z)=\Phi(z)M(\gamma).
\end{equation*}
\begin{Remark}
Thus one obtains a representation
$\gamma\mapsto M^{-1}(\gamma)\in \SL(2,\mathbb C)$ of the
fundamental group of the sphere with three punctures. This
representation is defined up to a conjugation, which is due to the
choice of
$a$ and $R_0$. We keep in mind this freedom and choose $\Phi(z)$ to be
the canonical solution $\Phi(z)=\Phi^{(0)}(z)$ in $z=0$.
\end{Remark}

Let $\gamma_0,\gamma_1,\gamma_\infty$ denote the usual set of
generators of the fundamental group $\pi_1(\widehat{\mathbb
C}\setminus \{0,1,\infty\})$, i.~e.\ positively oriented loops
around the points $0,1,\infty$. Denote by
\begin{equation*}
M_\nu :=M(\gamma_\nu),
\quad \nu=0,1,\infty,
\end{equation*}
the cor\-re\-spon\-ding mo\-no\-dro\-my ma\-tri\-ces
ge\-ne\-ra\-ting the mo\-no\-dro\-my group.\linebreak They satisfy
the cyclic relation
\begin{equation}                                            \label{cyclic}
M_\infty M_1 M_0=I.
\end{equation}
The canonical solutions differ by the \emph{connection matrices}
$E_\nu$
\begin{equation*}
\Phi^{(0)}(z)=\Phi^{(\nu)}(z)E_\nu,\qquad \nu=0,1,\infty.
\end{equation*}
By definition, $E_0=I$. Formulas for other two connection matrices
are more complicated and are proved in Appendix~\ref{proofs}.
\begin{Lemma}                                           \label{E_matrix}
The connection matrices are as follows:
\begin{align}
\label{E_1}
E_1&=
\begin{pmatrix}
-\frac{2\alpha+1}{\beta+\tau}
\frac{\Gamma(c)\Gamma(c-a-b)}{\Gamma(c-a)\Gamma(c-b)}&
\frac{2\alpha-1}{\gamma}
\frac{\Gamma(-c)\Gamma(c-a-b)}{\Gamma(-a)\Gamma(-b)}\\
&\\
-\frac{2\alpha+1}{\delta}
\frac{\Gamma(c)\Gamma(a+b-c)}{\Gamma(a)\Gamma(b)} e^{2\tau\pi i}&
-\frac{2\alpha-1}{\beta+\tau}
\frac{\Gamma(-c)\Gamma(a+b-c)}{\Gamma(a-c)\Gamma(b-c)} e^{2\tau\pi i}
\end{pmatrix},\\
\notag\\
E_\infty&=
\begin{pmatrix}
\frac{2\alpha+1}{\delta}\frac{a(\beta+\tau)}{\gamma(a-c)}
\frac{\Gamma(c)\Gamma(b-a)}{\Gamma(b)\Gamma(c-a)} e^{a\pi i}&
\frac{a(\beta+\tau)}{\gamma(a-c)}
\frac{\Gamma(2-c)\Gamma(b-a)}{\Gamma(b-c+1)\Gamma(1-a)} e^{(a-c)\pi i}\\
&\\
-\frac{2\alpha+1}{\delta}
\frac{\Gamma(c)\Gamma(a-b)}{\Gamma(a)\Gamma(c-b)} e^{b\pi i}&
-\frac{\Gamma(2-c)\Gamma(a-b)}{\Gamma(a-c+1)\Gamma(1-b)}
e^{(b-c)\pi i}
\end{pmatrix},
\end{align}
where the coefficients are as in Theorem~\ref{th_Fuchs}.
\end{Lemma}

The definition \eqref{can_sol} of the canonical solutions imply
for the monodromy matrices of the Fuchsian system
\eqref{sys_Fuchs}:
\begin{equation*}
M_\nu=E_\nu^{-1}e^{2\pi i\Lambda_\nu}E_\nu,
\qquad \nu=0,1,\infty.
\end{equation*}
Substituting formula \eqref{E_1}, and taking into account the
cyclic relation \eqref{cyclic}, and that the gauge transformation
\eqref{gauge} changes the sign of the monodromy matrix
at $z=1$, we arrive at the following theorem.
\begin{Theorem}
  The monodromy matrices of the solution $\Psi=D\Phi^{(0)}$ of the
  differential equation \eqref{WR_PQ} with $\spinorP$, $\spinorQ$ as in
  \eqref{W_data} are as follows:
  \begin{equation}           \label{mathcal_M}
    \begin{aligned}
      \mathcal M_0&=
      \begin{pmatrix}
        e^{2\pi i\alpha}&0\\
        0&e^{-2\pi i\alpha}
      \end{pmatrix},\qquad
      \mathcal M_\infty=\mathcal M_0^{-1}\mathcal M_1^{-1},\\
      \mathcal M_1&=
      \begin{pmatrix}
        e^{2\pi i\tau}-2i\frac{\sin\pi a\sin\pi b}{\sin\pi c} &\frac{2\pi
          i}{\gamma}\frac{2\alpha-1}{2\alpha+1}
        \frac{\Gamma^2(-c)}
        {\Gamma(-a)\Gamma(-b)\Gamma(a-c)\Gamma(b-c)}
        \\
        &\\
        \frac{2\pi i}{\delta}\frac{2\alpha+1}{2\alpha-1}
        \frac{\Gamma^2(c)}
        {\Gamma(a)\Gamma(b)\Gamma(c-a)\Gamma(c-b)} &e^{2\pi i\tau}+2i
        \frac{\sin\pi(c-a)\sin\pi(c-b)}{\sin\pi c}
      \end{pmatrix}.
    \end{aligned}
  \end{equation}
\end{Theorem}

\section{Trinoids. Moduli}

The immersion formula
\begin{equation}                                        \label{immersion}
  F=\Psi \Psi^*
\end{equation}
with $\Psi$ satisfying the trinoid equation describes a trinoid if and only
if the monodromy group of $\Psi$ is unitary. The monodromy group of the
equation is defined up to a conjugation. We call the monodromy group
\eqref{mathcal_M} \emph{unitarizable} if there exists $R\in \SL(2,{\mathbb
  C})$ such that all the matrices
\begin{equation*}
R^{-1}\mathcal M_0 R,\quad R^{-1}\mathcal M_1 R,\quad
R^{-1}\mathcal M_\infty R.
\end{equation*}
are unitary. In this case the immersion formula \eqref{immersion}
with $\Psi$ given by
\begin{equation*}
\Psi(z)=D\Phi^{(0)}(z) R
\end{equation*}
describes a trinoid.
\begin{Theorem}     \label{t.unitarizable}
  The monodromy group of the
  diffe\-rential equa\-tion
  \eqref{WR_PQ}, with $\spinorP$, $\spinorQ$ as in \eqref{W_data} is
  unitarizable if and only if
\begin{description}
  \item[(i)] $\alpha,\,\tau,\,\rho\,\in\,\mathbb R$;
  \item[(ii)] $
\sin \pi a\sin \pi b\sin \pi(a-c)\sin \pi (b-c)<0$.
\end{description}
\end{Theorem}
\begin{proof}
The necessity of the condition (i) is obvious. Then $a,b,c\in
{\mathbb R}$ and formulas of Theorem~\ref{th_gauge} imply
$\beta,\,\gamma,\,\delta\in \mathbb R$. The matrix
$R^{-1}{\mathcal M}_\nu R$ is unitary if and only if
$$
R R^*= {\mathcal M}_\nu R R^* {\mathcal M}_\nu^*.
$$
For $\nu=0$ this implies
$R R^*=\diag (r^2,r^{-2})$ with $r\in{\mathbb R}$.
For $\nu=1$ we get
\begin{equation}                                        \label{r}
r^4
\frac{(2\alpha+1)^2\gamma}{(2\alpha-1)^2\delta}
\frac{\Gamma^2(c)\Gamma(-a)\Gamma(-b)\Gamma(a-c)\Gamma(b-c)}
{\Gamma^2(-c)\Gamma(a)\Gamma(b)\Gamma(c-a)\Gamma(c-b)}=1.
\end{equation}
Applying
$
\Gamma(x)\Gamma(1-x)=\dfrac{\pi}{\sin(\pi x)}
$
we get that the right hand side is positive, and thus a formula
for
$r$, if and only if condition
(ii) holds.
\end{proof}

Umehara and Yamada \cite{Umehara/Yamada:2000} classify CMC-1 trinoids
according the conical singularities of their metric (see also
\cite{Umehara/Yamada:1996}). A conformal metric $e^u\,dz\,d\overline{z}$ is
said to have a {\em conical singularity of order $\beta$} at $z=z_0$, if
\begin{equation*}
  u=2\beta\log|z-z_0| + O(z-z_0)\quad \text{for}\quad z\to z_0.
\end{equation*}
The metric of a CMC-1 trinoid has three conical singularities. Let their
degrees be $\beta_1$, $\beta_2$, and $\beta_3$, and let $B_j=\pi(\beta_j+1)$.
Umehara and Yamada derive the following condition for the $B_j$:
\begin{equation*}
\cos^2 B_1 + \cos^2 B_2 + \cos^2 B_3 + 2\cos B_1 \cos B_2 \cos B_3 < 1.
\end{equation*}
It turns out that this condition is equivalent to inequality
(ii) of theorem \ref{t.unitarizable}.  The crux is to show that $\beta_1 =
-2(\alpha+k_1)$, $\beta_2=-2(\tau+k_2)$, and $\beta_3=-2(\rho+k_3)$ with
$k_1, k_2, k_3 \in\mathbb{Z}$. From this, one obtains by a long but
elementary calculation
\begin{multline*}
  \cos^2 B_1 + \cos^2 B_2 + \cos^2 B_3 + 2\cos B_1 \cos B_2 \cos B_3 - 1\\
  = \sin \pi a\sin \pi b\sin \pi(a-c)\sin \pi (b-c).
\end{multline*}
Indeed, equation \eqref{Psi} expresses the unbranched $\Phi_1$ as the product
of $\Psi$ and $\Phi_0$. Since $\Psi z^{-\Lambda_0}$ is unbranched at $z=0$,
this implies that $z^{\Lambda_0}\Phi_0$ is also unbranched at $z=0$. From
this, one deduces that $z^{\alpha}\spinorp$ and $z^{\alpha}\spinorq$ are
meromorphic at $0$. With \eqref{metric&hopf}, one obtains $\beta_1 =
-2(\alpha+k_1)$. The analogous expressions for $\beta_2$ and $\beta_3$ follow
by symmetry.

We will derive a condition in terms of the parameters $p_0$, $p_1$,
$p_{\infty}$, $q_0$, $q_1$, $q_{\infty}$ of the Weierstrass data
\eqref{W_data}.  It is convenient to shift $c$'s in \eqref{c} by $1/4$,
\begin{align*}
d_0&=\frac{1}{4}+\oz+\zinf,\\
d_1&=\frac{1}{4}+\oz+\oinf,\\
d_\infty&=\frac{1}{4}+\zinf+\oinf.
\end{align*}
Introduce the fractional part $\{ x\}$ as a mapping
$$
\{\,\}:\,\mathbb R\to\,[-\tfrac{1}{2},\tfrac{1}{2}).
$$

\begin{Proposition}                               \label{c.unitaryModuli}
  The monodromy group of the diffe\-rential equa\-tion \eqref{WR_PQ}, with
  $\spinorP$, $\spinorQ$ as in \eqref{W_data} is uni\-ta\-ri\-zable if and
  only if $d_0,d_1,d_\infty \ge 0$ and
  $$(|\{\sqrt{d_0}\}|,|\{\sqrt{d_1}\}|,|\{\sqrt{d_\infty}\}|)\in{\mathcal
    D},$$
  where
  \begin{equation}                                            \label{moduli}
    \mathcal D=
    \left\{
      (\Delta_1,\Delta_2,\Delta_3)\in \mathbb R^3\;:
      \Delta_1,\Delta_2,\Delta_3\ge 0,\quad
      \begin{aligned}
        &\Delta_1+\Delta_2+\Delta_3> \tfrac{1}{2}\,,\\
        &\Delta_1+\Delta_2-\Delta_3< \tfrac{1}{2}\,,\\
        &\Delta_1+\Delta_3-\Delta_2< \tfrac{1}{2}\,,\\
        &\Delta_2+\Delta_3-\Delta_1< \tfrac{1}{2}\,.
      \end{aligned}
    \right\}
  \end{equation}
\end{Proposition}
\begin{proof}
  Condition (i) of Theorem~\ref{t.unitarizable} and
  $\alpha=\tfrac{1}{2}-\sqrt{d_0},\,\tau=\sqrt{d_1},\,\rho=\sqrt{d_\infty}$
  imply that $d_0, d_1, d_\infty$ are non-negative. We have
  \begin{equation*}
    a=\frac{1}{2}-\sqrt{d_0}+\sqrt{d_1}+\sqrt{d_\infty},\
    b=\frac{1}{2}-\sqrt{d_0}+\sqrt{d_1}-\sqrt{d_\infty},\
    c=1-2\sqrt{d_0}.
  \end{equation*}
  Further, using
  \begin{equation*}
    \begin{aligned}
      &\sin\pi a\,\sin\pi b\,\sin\pi(a-c)\,\sin\pi(b-c)=\\
      &\left(\cos\pi(a-b)-\cos\pi(a+b)\right)\left(\cos\pi(a-b)-\cos\pi(a+b-2c)\right)=\\
      &\left(\cos2\pi\sqrt{d_\infty}+\cos2\pi(\sqrt{d_0}-\sqrt{d_1})\right)
      \left(\cos2\pi\sqrt{d_\infty}+\cos2\pi(\sqrt{d_0}+\sqrt{d_1})\right)=\\
      &\cos\pi(|\{\sqrt{d_0}\}|+|\{\sqrt{d_1}\}|+|\{\sqrt{d_\infty}\}|)
      \cos\pi(|\{\sqrt{d_0}\}|+|\{\sqrt{d_1}\}|-|\{\sqrt{d_\infty}\}|)\times\\
      &\cos\pi(|\{\sqrt{d_0}\}|-|\{\sqrt{d_1}\}|+|\{\sqrt{d_\infty}\}|)
      \cos\pi(-|\{\sqrt{d_0}\}|+|\{\sqrt{d_1}\}|+|\{\sqrt{d_\infty}\}|)
    \end{aligned}
  \end{equation*}
  we transform condition (ii) of Theorem~\ref{t.unitarizable} to
  \hbox{$(|\{\sqrt{d_0}\}|,|\{\sqrt{d_1}\}|,|\{\sqrt{d_\infty}\}|)\in{\mathcal D}$}.
\end{proof}

There is an elementary derivation of the description \eqref{moduli} of the
moduli space, which does not involve hypergeometric functions.
Indeed, the local data provide us with the local monodromies (\ref{loc_mon}),
i.~e.\ with the eigenvalues
$$
e^{\pm 2\pi \alpha_0},\qquad e^{\pm 2\pi \alpha_1},
\qquad e^{\pm 2\pi \alpha_\infty}
$$
of the monodromy matrices ${\mathcal M}_0, {\mathcal M}_1$ and
${\mathcal M}_\infty$.
We have ${\mathcal M}_0, {\mathcal M}_1, {\mathcal M}_\infty\in \SL(2,
{\mathbb C})$ with ${\mathcal M}_\infty {\mathcal M}_1 {\mathcal M}_0=I$,
and the problem is to characterize the unitarizable monodromies. This problem
is equivalent to the following one: What is the necessary and sufficient
condition for the existence of a gauge $G\in \SL(2, {\mathbb C})$ such that
all the matrices $G{\mathcal M}_0 G^{-1}, G{\mathcal M}_1 G^{-1}, G{\mathcal
  M}_\infty G^{-1}$ belong to $\SU(2)$?

First, let us normalize the eigenvalues as follows:
\begin{equation*}
0<\alpha_0,\alpha_1,\alpha_\infty<\frac{1}{2}.
\end{equation*}
Note that the case of half-integer coefficients is excluded
\eqref{generic}. Without loss of generality one can assume
${\mathcal M}_0$ to be diagonal
$$
{\mathcal M}_0=
\begin{pmatrix}
e^{2\pi i\alpha_0}&0
\\
0&e^{-2\pi i\alpha_0}
\end{pmatrix}.
$$
Further, by an appropriate diagonal gauge transformation let us
normalize the sum of the off-diagonal terms of ${\mathcal M}_1$ to
vanish:
$$
{\mathcal M}_1=
\begin{pmatrix}
u&v
\\
-v&w
\end{pmatrix},
\qquad
{\mathcal M}_\infty^{-1}={\mathcal M}_1 {\mathcal M}_0=
\begin{pmatrix}
u e^{2\pi i\alpha_0}& v e^{-2\pi i\alpha_0}
\\
-v e^{2\pi i\alpha_0}&we^{-2\pi i\alpha_0}
\end{pmatrix}.
$$
Now, ${\mathcal M}_1,{\mathcal M}_\infty\in \SU(2)$
if and only if
\begin{equation*}                        \label{unitarity_uvw}
\begin{split}
uw+v^2&=1,\\
u+w&=2\cos 2\pi\alpha_1,\\
u e^{2\pi i \alpha_0}+w e^{-2\pi i \alpha_0}&=2\cos
2\pi\alpha_\infty
\end{split}
\end{equation*}
with real $v\in {\mathbb R}$. The last two equations are
equivalent to $w=\bar{u}$ and
\begin{eqnarray*}
{\rm Re}\,u &=&\cos 2\pi\alpha_1,\\
\sin 2\pi\alpha_0 {\rm Im}\,u&=&\cos 2\pi\alpha_0 \cos 2\pi\alpha_1-\cos 2\pi\alpha_\infty.
\end{eqnarray*}
There exists real $v$ in the first equation of
(\ref{unitarity_uvw}) if and only if
$$
|{\rm Im}\,u|<\sin 2\pi\alpha_1.
$$
Substituting the formula for ${\rm Im}\,u$ we obtain the system
\begin{eqnarray*}
\cos 2\pi(\alpha_0+\alpha_1)&<& \cos 2\pi\alpha_\infty,\\
\cos 2\pi(\alpha_0-\alpha_1)&>& \cos 2\pi\alpha_\infty.
\end{eqnarray*}
With the chosen normalization these two inequalities are
equivalent to
\begin{eqnarray*}
1-\alpha_\infty>\alpha_0+\alpha_1>\alpha_\infty, \\
\alpha_\infty>\alpha_0-\alpha_1>-\alpha_\infty
\end{eqnarray*}
 respectively. Finally we get the
following conditions
\begin{equation}                       \label{unitarity_alpha}
\begin{aligned}
\alpha_0+\alpha_1+\alpha_\infty &<&1,\\
\alpha_0+\alpha_1-\alpha_\infty &>&0,\\
\alpha_0-\alpha_1+\alpha_\infty &>&0,\\
-\alpha_0+\alpha_1+\alpha_\infty &>&0.\\
\end{aligned}
\end{equation}
In our notations (\ref{loc_mon}) we have
$$
\alpha_i\equiv\pm (\frac{1}{2}+\sqrt{d_i})\qquad ({\rm mod}\ {\mathbb Z}).
$$
The representative in the interval $(0,\frac{1}{2})$ is
$$
\alpha_i= \frac{1}{2}-|\{\sqrt{d_i})\}|.
$$
Finally, written in terms of $|\{\sqrt{d_i})\}|$ conditions
(\ref{unitarity_alpha}) coincide with (\ref{moduli}).

Note that condition (\ref{moduli}) can be derived from a result of
Biswas \cite{Biswas}. He found the necessary and sufficient
condition for the existence of a flat irreducible $U(2)$
connection on a punctured sphere such that the local monodromies
around any puncture is in the preassigned conjugacy class. Since
for three punctures the conjugacy classes data determine the
monodromy,
Biswas' condition characterizes the unitarizable monodromy groups
of trinoids.

Finally, CMC-1 trinoids are constructed as follows: Take Weierstrass data
$p_0,p_1,p_\infty, q_0,q_1,q_\infty$ satisfying the conditions of
proposition~\ref{c.unitaryModuli} and apply the immersion formula
\eqref{immersion} with $\Psi$ given by
\begin{equation*}
\begin{aligned}
\Psi(z)=&D\Phi^{(0)}(z)R\\
       =&D\Phi^{(1)}(z)E_1 R\\
       =&D\Phi^{(\infty)}(z)E_\infty R,
\end{aligned}
\end{equation*}
choosing the representation converging in the corresponding
parameter domain. Here, $D$ is the gauge matrix \eqref{gauge},
$\Phi^{(\nu)}(z)$  the canonical solutions \eqref{Phi_0},
\eqref{Phi_1}, \eqref{Phi_inf}, $E_\nu$ the connection matrices
and $R=\diag (r,r^{-1})$ with $r$ from \eqref{r}.


\begin{figure}[ht]
\centering
\mbox{\subfigure[$d_0<D_0$]{\epsfig{figure=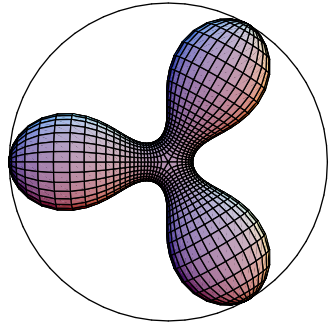,width=.40\textwidth}}\quad
    \subfigure[$d_0=0,2332 \approx D_0$]{\epsfig{figure=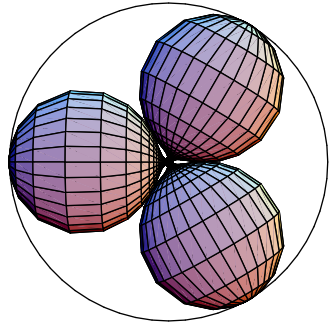,width=.40\textwidth}}
}
\mbox{\subfigure[$d_0=0,2400$]{\epsfig{figure=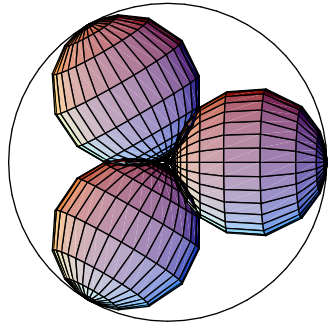,width=.40\textwidth}}\quad
    \subfigure[$d_0>D_0$]{\epsfig{figure=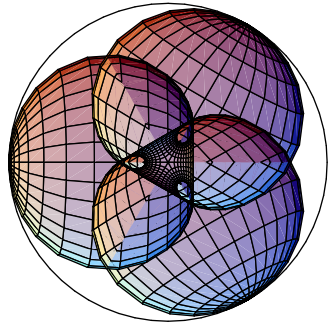,width=.40\textwidth}}
}
\caption{Symmetric trinoids \label{f.symmetric}}
\end{figure}

The CMC-1 trinoids build a three-parameter family. Indeed, let us fix the
points on the absolute applying isometries of $\mathbb H^3$. We
fix the images of the ends
$z_j, j=0,1,\infty$ at the points
$(-\tfrac{1}{2},0,-\tfrac{\sqrt{3}}{2})$, $(1,0,0)$ and
$(-\tfrac{1}{2},0,\tfrac{\sqrt{3}}{2})$ respectively.
This implies the following relations between the parameters
$p_0,p_1,p_\infty, q_0,q_1,q_\infty$:
\begin{equation}
\begin{aligned}
p_0 =&(2-\sqrt{3})q_0,\\
p_1 =&-q_1,\\
p_\infty=&(2+\sqrt{3})q_\infty.
\end{aligned}
\label{param.count}
\end{equation}
Two embedded examples are shown in Fig.~\ref{f.non-symmetric} in
the introduction.

The condition $d_0=d_1=d_\infty$
characterizes the symmetric trinoids. They build a one-parameter
family characterized by the parameter $d_0$. There is
$D_0<\tfrac{1}{4}$ such that all symmetric trinoids with $d_0<D_0$ are
embedded and all trinoids with $d_0>D_0$ are not embedded, see
Fig.~\ref{f.symmetric}.

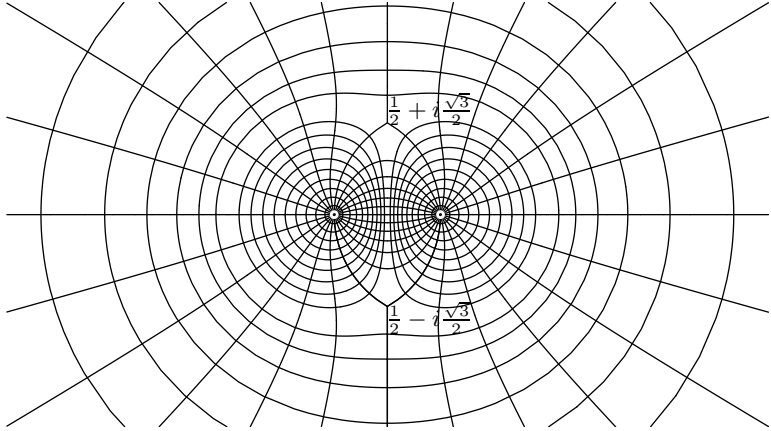
\begin{figure}[ht]
\begin{center}
\input{param.pstex_t}
\end{center}
\caption{Parameter lines in the $z$-plane
\label{f.parametr}}
\end{figure}

Figs.~\ref{f.non-symmetric} and \ref{f.symmetric} were produced using the
software \emph{Mathematica}, which provides an implementation of the
hypergeometric function.
The parameterization in these figures was chosen
to show the umbilic points at the centers of both sides of the symmetric
trinoids. We split the complex $z$-plane in three domains associated to the
ends and use the following parameterization for these domains
\begin{equation*}
z=
  \begin{cases}
    z(w_0),&\text{for end $z=0$},\\
    z(w_1),&\text{for end $z=1$},\\
    z(w_\infty),&\text{for end $z=\infty$}
  \end{cases}\quad
  z(w_j)=\frac{z_1-z_\infty}{z_1-z_0}\frac{w_j-z_0}{w_j-z_\infty},
\end{equation*}
where
$w_j=(-\frac{\tilde w+1}{\tilde w-1})^{\tfrac{2}{3}}z_j$,
$z_0=e^{i \tfrac{\pi}{6}}$,
$z_1=e^{i \tfrac{5 \pi}{6}}$,
$z_\infty=e^{i \tfrac{3 \pi}{2}}$,
$|\tilde w|\le 1$.
The corresponding parameter lines $w\in \mathbb R$, $w\in i\mathbb
R$ in the $z$-plane are shown in Fig.~\ref{f.parametr}.

The \emph{Mathematica} notebook as well as additional images of
trinoids can be found from the URL {\tt
http://www-sfb288.math.tu-berlin.de/\~\,bobenko}
\vspace{5mm}

{\bf Acknowledgment} We are grateful to Alexander Its for a useful
correspondence and for showing his formulas on integration of
Fuchsian systems with three singularities in hypergeometric
functions prior to publication of his book. This research was
financially supported by Deutsche Forschungsgemeinschaft (SFB 288
"Differential Geometry and Quantum Physics").

\appendix
\section{Basic facts about the hypergeometric function}

We present some facts about the hypergeometric function used
in the proofs of Theorem~\ref{th_Fuchs} and Lemma~\ref{E_matrix}.

The function represented by the infinite series
$\overset{^\infty}{\underset{n=0}{\sum}} \frac{(a)_n\,(b)_n}{(c)_n}\,
\frac{z^n}{n!}$ within its circle of convergence and its analytic
continuation is called \emph{the hypergeometric function} $\F
(a,b;c;z)$. The symbol $(a)_n$ is defined as
\begin{equation*}
(a)_n\,=\,a(a+1)(a+2)\dotsb(a+n-1)\,=\,\frac{\Gamma(a+n)}{\Gamma(a)}\,.
\end{equation*}

Thus
\begin{equation*}
\F(a,b;c;z)\,=\,\frac{\Gamma(c)}{\Gamma(a)\,\Gamma(b)}\,
\sum_{n=0}^{\infty}\frac{\Gamma(a+n)\,\Gamma(b+n)}{\Gamma(c+n)}\frac{z^n}{n!}\,.
\end{equation*}

For later use we give some formulas for hy\-per\-geo\-met\-ric functions
\cite{Magnus/Ober/Soni, Bateman/Erdelyi:Vol_1, Whittaker/Watson}.

\paragraph{Differentiation formula.}
\begin{equation*}
\frac{d^n}{dz^n}\,\F(a,b;c;z)\,=\,\frac{(a)_n\,(b)_n}{(c)_n}\,
\F(a+n,b+n;c+n;z)\,.
\end{equation*}

\paragraph{Gauss' contiguous relations.}
\noindent The six functions
$$\F(a\pm 1,b;c;z),\quad \F(a,b\pm 1;c;z),\quad \F(a,b;c\pm 1;z)$$
are called
\emph{contiguous} to $\F(a,b;c;z)$.  A relation between $\F(a,b;c;z)$ and any
two contiguous functions is called \emph{a contiguous relation}. By these
relations, one can expresses the function $\F(a+l,b+m;c+n;z)$ with
$l,m,n\in\mathbb{Z}$, $c+n\ne 0,\,-1,\,-2,\,\dotsc$ as a linear combination
of $\F(a,b;c;z)$ and one of its contiguous functions.  The coefficients are
rational functions of $a,\,b,\,c,\,z$. For example, one has the following
formulas:
\begin{multline}\label{ReduceGauss1}
\F(a+1,b+1;c+1;z)\,=\,\frac{1}{a\,b\,(1-z)}\, \\
\left[c\,(a+b-c)\;\F(a,b;c;z)+(c-a)\,(c-b)\;\F(a,b;c+1;z)\right]\,,
\end{multline}
\begin{multline}\label{ReduceGauss2}
\F(a+1,b+1;c+2;z)\,=\,\\
\frac{c\,(c+1)}{a\,b\,z}\,
\left[\F(a,b;c;z)-\F(a,b;c+1;z)\right]\,.
\end{multline}

\paragraph{The connection between hypergeometric functions of $z$ and of
$1-z$.}
For $|\arg(1-z)|<\pi$ and $c-a-b\not\in\{0,\pm 1,\pm 2,\dotsc\}$,
\begin{multline}\label{connect0_1}
    \F(a,b;c;z) = \\
    \frac{\Gamma(c)\Gamma(c-a-b)}{\Gamma(c-a)\Gamma(c-b)}
    \F(a,b;a+b-c+1;1-z)+\\
    (1-z)^{c-a-b}\frac{\Gamma(c)\Gamma(a+b-c)}{\Gamma(a)\Gamma(b)}
    \F(c-a,c-b;c-a-b+1;1-z).
\end{multline}

\paragraph{The hypergeometric differential equation}
\begin{equation}
\label{HyperGeomDE}
z\,(1-z)\,\frac{d^2w}{dz^2}+[c-(a+b+1)\,z]\,\frac{dw}{dz}-a\,b\,w\,=\,0
\end{equation}
has three regular singular points $z=0,\,1,\,\infty$. The pairs of
characteristic exponents at these points are
\begin{alignat*}{3}
\rho_0\,&=\,0\,,&\quad \rho_1\,&=\,0\,,&\quad
\rho_{\infty}\,&=\,a\,,\\
\rho_0'\,&=\,1-c\,,&\quad \rho_1'\,&=\,c-a-b\,,&\quad
\rho_{\infty}'\,&=\,b
\end{alignat*}
respectively. The hypergeometric function $\F(a,b;c;z)$ is a
solution of the hypergeometric differential equation which
is unbranched at $z=0$.

\begin{Proposition}[\cite{Klein}]\label{SolSys}
The fundamental system of linearly independent solutions of
hypergeometric differential equation \eqref{HyperGeomDE} at the
singular points $z=0,1,\infty$ is given by
\begin{align*}
w_1^{(0)}(z)\,&=\,\F(a,b;c;z)\,,\\
w_2^{(0)}(z)\,&=\,z^{1-c}\,\F(a-c+1,b-c+1;2-c;z)\,,\\
w_1^{(1)}(z)\,&=\,\F(a,b;a+b-c+1;1-z)\,,\\
w_2^{(1)}(z)\,&=\,(1-z)^{c-a-b}\,\F(c-b,c-a;c-a-b+1;1-z)\,,\\
w_1^{(\infty)}(z)\,&=\,z^{-a}\,\F(a,a-c+1;a-b+1;\frac{1}{z})\,,\\
w_2^{(\infty)}(z)\,&=\,z^{-b}\,\F(b,b-c+1;b-a+1;\frac{1}{z})\,.
\end{align*}
\end{Proposition}

\paragraph{Riemann's differential equation.}
The hypergeometric differential equation is a special case of
Riemann's differential equation
\begin{multline*}
\frac{d^2w}{dz^2}+
\left[
\frac{1-\rho_a-\rho_a'}{z-a}+
\frac{1-\rho_b-\rho_b'}{z-b}+
\frac{1-\rho_c-\rho_c'}{z-c}
\right]\,\frac{dw}{dz}\\
+\left[
\frac{\rho_a\,\rho_a'\,(a-b)\,(a-c)}{z-a}+
\frac{\rho_b\,\rho_b'\,(b-c)\,(b-a)}{z-b}+
\frac{\rho_c\,\rho_c'\,(c-a)\,(c-b)}{z-c}
\right]\\
\times \frac{w}{(z-a)\,(z-b)\,(z-c)}\,=\,0\,.
\end{multline*}
The characteristic exponents $\rho_a$, ${\rho_a}'$; ${\rho_b}$, ${\rho_b}'$;
$\rho_c$, ${\rho_c}'$ must satisfy the additional relation
\begin{equation*}
\sum_{j=0,1,\infty}\,(\rho_j+{\rho_j}') = 1.
\end{equation*}
The following symbol is used for Riemann's differential equation:
\begin{equation*}
P\left\{
\begin{matrix}
a&b&c&\\
\rho_a&\rho_b&\rho_c&;z\\
\rho_a'&\rho_b'&\rho_c'&
\end{matrix}
\right\}\,.
\end{equation*}
It is also used to denote the set of solutions of the equation and called
\emph{Riemann P-function}.

In particular, the hypergeometric  differential equation
\eqref{HyperGeomDE} is
\begin{equation*}
P\left\{
\begin{matrix}
0&\infty&1&\\
0&a&0&;z\\
1-c&b&c-a-b&
\end{matrix}
\right\}\,.
\end{equation*}
The generalized hypergeometric  differential equation
\begin{multline}
\label{GenHyperGeomDE}
\frac{d^2w}{dz^2}+
\left[
\frac{1-\rho_0-\rho_0'}{z}+
\frac{1-\rho_1-\rho_1'}{z-1}
\right]\,\frac{dw}{dz}\\
+\left[
\frac{-\rho_0\,\rho_0'}{z}+
\frac{\rho_1\,\rho_1'}{z-1}+
\rho_{\infty}\,\rho_{\infty}'
\right]\,\frac{w}{z\,(z-1)}\,=\,0
\end{multline}
is represented by
\begin{equation*}
P\left\{
\begin{matrix}
0&\infty&1&\\
\rho_0&\rho_{\infty}&\rho_1&;z\\
\rho_0'&\rho_{\infty}'&\rho_1'&
\end{matrix}
\right\}\,.
\end{equation*}
The following two transformations for\-mu\-las are valid for
Riemann's
$P$-function

\begin{multline*}
1.\quad P\left\{
\begin{matrix}
0&\infty&1&\\
\rho_0&\rho_{\infty}&\rho_1&;z\\
\rho_0'&\rho_{\infty}'&\rho_1'&
\end{matrix}
\right\}\\
=\,z^{-k}\,(z-1)^{-l}\,P\left\{
\begin{matrix}
0&\infty&1&\\
\rho_0+k&\rho_{\infty}-k-l&\rho_1+l&;z\\
\rho_0'+k&\rho_{\infty}'-k-l&\rho_1'+l&
\end{matrix}
\right\}\,;\\
\shoveleft{2a.\quad
P\left\{
\begin{matrix}
0&\infty&1&\\
\rho_0&\rho_{\infty}&\rho_1&;z\\
\rho_0'&\rho_{\infty}'&\rho_1'&
\end{matrix}
\right\}\,
=\, P\left\{
\begin{matrix}
1&\infty&0&\\
\rho_0&\rho_{\infty}&\rho_1&;1-z\\
\rho_0'&\rho_{\infty}'&\rho_1'&
\end{matrix}
\right\}\,;}\\
\shoveleft{2b.\quad
P\left\{
\begin{matrix}
0&\infty&1&\\
\rho_0&\rho_{\infty}&\rho_1&;z\\
\rho_0'&\rho_{\infty}'&\rho_1'&
\end{matrix}
\right\}\,
=\, P\left\{
\begin{matrix}
\infty&0&1&\\
\rho_0&\rho_{\infty}&\rho_1&;\frac{1}{z}\\
\rho_0'&\rho_{\infty}'&\rho_1'&
\end{matrix}
\right\}\,.
\hfill}
\end{multline*}
\begin{Proposition}\label{GHDESystem}
The fundamental system of linear independent solutions of
generalized hypergeometric equation \eqref{GenHyperGeomDE} at the
singular points $z=0,\,1,\,\infty$ is given by
\begin{align*}
w_1^{(0)}(z)\,&=\,z^{\rho_0}\,(z-1)^{\rho_1}\,\F(a,b;c;z)\,,\\
w_2^{(0)}(z)\,&=\,z^{\rho_0'}\,(z-1)^{\rho_1}\,\F(a-c+1,b-c+1;2-c;z)\,,\\
w_1^{(1)}(z)\,&=\,z^{\rho_0}\,(z-1)^{\rho_1}\,\F(a,b;a+b-c+1;1-z)\,,\\
w_2^{(1)}(z)\,&=\,z^{\rho_0}\,(z-1)^{\rho_1'}\,\F(c-b,c-a;c-a-b+1;1-z)\,,\\
w_1^{(\infty)}(z)\,&=\, z^{-\rho_1-\rho_\infty}\,(z-1)^{\rho_1}\,
\F(a,a-c+1;a-b+1;\frac{1}{z})\,,\\
w_2^{(\infty)}(z)\,&=\, z^{-\rho_1-\rho_\infty'}\,(z-1)^{\rho_1}\,
\F(b,b-c+1;b-a+1;\frac{1}{z})\,,
\end{align*}
where
\begin{equation*}
a\,=\,\rho_0+\rho_1+\rho_\infty\,,\quad
b\,=\,\rho_0+\rho_1+\rho_\infty'\,,\quad c\,=\,1+\rho_0-\rho_0'\,.
\end{equation*}
\end{Proposition}

\section{The proofs of Theorem~\ref{th_Fuchs} and Lemma~\ref{E_matrix}}\label{proofs}
\begin{proof}[Proof of Theorem~\ref{th_Fuchs}]
Reduce the Fuchsian system \eqref{sys_Fuchs}
\begin{align*}
\begin{pmatrix}
\phi_{11}'&\phi_{12}'\\
\phi_{21}'&\phi_{22}'
\end{pmatrix}\,&=\,
\begin{pmatrix}
\frac{\alpha}{z}+\frac{\beta}{z-1}&\frac{\gamma}{z-1}\\
\frac{\epsilon^2\,\delta}{z-1}&-\frac{\alpha}{z}-\frac{\beta}{z-1}
\end{pmatrix}\cdot
\begin{pmatrix}
\phi_{11}&\phi_{12}\\
\phi_{21}&\phi_{22}
\end{pmatrix}\,\\
&=\,
\begin{pmatrix}
a_{11}&a_{12}\\
a_{21}&a_{22}
\end{pmatrix}\cdot
\begin{pmatrix}
\phi_{11}&\phi_{12}\\
\phi_{21}&\phi_{22}
\end{pmatrix}
\end{align*}
to the following system of second-order differential equations.
\begin{equation*}
\left\{
\begin{minipage}[c]{0.97\textwidth}
\begin{align}
\phi_{1j}'
&=a_{11}\phi_{1j}+a_{12}\phi_{2j},\quad j=1,2 \label{first1j}\\
\phi_{2j}'
&=a_{21}\phi_{1j}+a_{22}\phi_{2j},\quad j=1,2 \label{first2j}\\
\phi_{1j}''
&=\phi_{1j}'\left(a_{11}+a_{22}+\tfrac{a_{12}'}{a_{12}}\right)+
\phi_{1j}\left(a_{12}a_{21}-a_{11}a_{22}
+\tfrac{a_{11}'a_{12}-a_{11}a_{12}'}{a_{12}}\right)\notag\\
&=\phi_{1j}'\left(-\tfrac{1}{z-1}\right)
+\phi_{1j}\left(\tfrac{\alpha^2-\alpha}{z^2}+
\tfrac{\beta^2+\gamma\delta}{(z-1)^2}+
\tfrac{\alpha(2\beta+1)}{z(z-1)}\right),\quad j=1,2 \label{second1j} \\
\phi_{2j}''
&=\phi_{2j}'\left(a_{11}+a_{22}+\tfrac{a_{21}'}{a_{21}}\right)+
\phi_{2j}\left(a_{12}a_{21}-a_{11}a_{22}
+\tfrac{a_{22}'a_{21}-a_{22}a_{21}'}{a_{21}}\right)\notag\\
&=\phi_{2j}'\left(-\tfrac{1}{z-1}\right)
+\phi_{2j}\left(\tfrac{\alpha^2+\alpha}{z^2}+
\tfrac{\beta^2+\gamma\delta}{(z-1)^2}+
\tfrac{\alpha(2\beta-1)}{z(z-1)}\right),\quad j=1,2.\label{second2j}
\end{align}
\end{minipage}
\right.\notag
\end{equation*}
Equations \eqref{second1j}, \eqref{second2j} are the generalized
hypergeometric differential equations with the characteristic
exponents
\begin{alignat*}{3}
\rho_0\,&=\,\alpha\,,&\quad
\rho_1\,&=\,\sqrt{\beta^2+\gamma\,\delta}\,,&\quad
\rho_{\infty}\,&=\,\sqrt{(\alpha+\beta)^2+\gamma\,\delta}\,,\\
\rho_0'\,&=\,1-\alpha\,,&\quad
\rho_1'\,&=\,-\sqrt{\beta^2+\gamma\,\delta}\,,&\quad
\rho_{\infty}'\,&=\,-\sqrt{(\alpha+\beta)^2+\gamma\,\delta}\\
\intertext{and}
\widetilde \rho_0\,&=\,1+\alpha\,,&\quad
\widetilde \rho_1\,&=\,\sqrt{\beta^2+\gamma\,\delta}\,,&\quad
\widetilde \rho_{\infty}\,&=\,\sqrt{(\alpha+\beta)^2+\gamma\,\delta}\,,\\
\widetilde \rho_0'\,&=\,-\alpha\,,&\quad
\widetilde \rho_1'\,&=\,-\sqrt{\beta^2+\gamma\,\delta}\,,&\quad
\widetilde \rho_{\infty}'\,&=\,-\sqrt{(\alpha+\beta)^2+\gamma\,\delta}
\end{alignat*}
respectively. Chose the ansatz
\begin{equation*}
\Phi^{(0)}(z)\,=\,
\begin{pmatrix}
k_{11}\,w_1^{(0)}(z)&w_2^{(0)}(z)\\
\widetilde w_1^{(0)}(z)&k_{22}\,\widetilde w_2^{(0)}(z)
\end{pmatrix}\,,
\end{equation*}
for a solution of the Fuchsian system at $z=0$. Here,
$w_1^{(0)}(z),\,w_2^{(0)}(z)$ and
$\widetilde w_1^{(0)}(z),\,\widetilde w_2^{(0)}(z)$ are
linearly independent solutions of equations \eqref{second1j} and
\eqref{second2j}, respectively. Due to
Proposition~\ref{GHDESystem}, the function
$w_1^{(0)}(z)$, $w_2^{(0)}(z)$, $\widetilde w_1^{(0)}(z)$,
$\widetilde w_2^{(0)}(z)$ can be chosen as follows:
\begin{align*}
w_1^{(0)}(z)
&=\,z^{\alpha}\,(z-1)^{\tau}\,\F(a,b;c;z)\,,\\
w_2^{(0)}(z)
&=\,z^{1-\alpha}\,(z-1)^{\tau}\,\F(a-c+1,b-c+1;2-c;z)\,,\\
\widetilde w_1^{(0)}(z)
&=\,z^{1+\alpha}\,(z-1)^{\tau}\,\F(a+1,b+1;c+2;z)\,,\\
\widetilde w_2^{(0)}(z)
&=\,z^{-\alpha}\,(z-1)^{\tau}\,\F(a-c,b-c;-c;z)\,,
\end{align*}
where
\begin{align*}
a\,&=\,\rho_0+\rho_1+\rho_\infty\,=\,\alpha+\tau+\rho\,,\\
b\,&=\,\rho_0+\rho_1+\rho_\infty'\,=\,\alpha+\tau-\rho\,,\\
c\,&=\,1+\rho_0-\rho_0'\,=\,2\alpha\,,\\
\widetilde a\,&=\,\widetilde \rho_0+\widetilde \rho_1+\widetilde \rho_\infty\,
=\,1+\alpha+\tau+\rho\,=\,a+1,\\
\widetilde b\,&=\,\widetilde \rho_0+\widetilde \rho_1+\widetilde \rho_\infty'\,
=\,1+\alpha+\tau-\rho\,=\,b+1,\\
\widetilde c\,&=\,1+\widetilde \rho_0-\widetilde \rho_0'\,
=\,2+2\alpha\,=\,c+2\,.
\end{align*}
The coefficients $k_{11},\,k_{22}$ follow from the conditions
\eqref{first1j}, \eqref{first2j}:
\begin{multline*}
k_{11}\,\frac{d}{dz}w_1^{(0)}(z)-a_{11}\,k_{11}\,w_1^{(0)}(z)
-a_{12}\,k_{21}\,\widetilde w_1^{(0)}(z)\,=\\
\shoveleft{z^\alpha\,(z-1)^{\tau-1}\,
\biggl[-\gamma\,z\,\F(a+1,b+1;c+2,z)\biggr.}\\
\biggl.+k_{11}\,\biggl((\tau-\beta)\,\F(a,b;c,z)
+\frac{a\,b}{c}\,(z-1)\,\F(a+1,b+1;c+1,z)\biggr)\biggr]
\overset{\eqref{ReduceGauss1}, \eqref{ReduceGauss2}}{=}\,\\
\shoveleft{z^\alpha\,(z-1)^{\tau-1}\,
\biggl[-\gamma\,\frac{c\,(c+1)}{a\,b}\,
\biggl(\F(a,b;c,z)-\F(a,b;c+1,z)\biggr)\biggr.}\\
\biggl.+k_{11}\,\biggl((\tau-\beta+a+b-c)\,\F(a,b;c,z)
-\frac{(c-a)\,(c-b)}{c}\,\F(a,b;c+1;z)\biggr)\biggr]\,=\\
\shoveleft{z^\alpha\,(z-1)^{\tau-1}\,
\biggl[\frac{\gamma\,(2\alpha+1)}{\tau-\beta}\,
\biggl(\F(a,b;c+1;z)-\F(a,b;c,z)\biggr)\biggr.}\\
\biggl.+k_{11}\,(\tau+\beta)\,
\biggl(\F(a,b;c+1;z)-\F(a,b;c,z)\biggr)\biggr]\,=\,0.
\end{multline*}
Hence, $k_{11}\,=\,-\frac{2\alpha+1}{\delta}$. The formula for
$k_{22}$ is obtained analogously.

From Proposition~\ref{GHDESystem} we know the system of two linear
independent solutions for the equations \eqref{second1j},
\eqref{second2j} both in the neighborhood of $z=1$ and of
$z=\infty$. In the same way as for the canonical solution
$\Phi^{(0)}(z)$ we prove the formulas for  $\Phi^{(1)}(z)$ and
$\Phi^{(\infty)}(z)$.
\end{proof}

\begin{proof}[Proof of Lemma~\ref{E_matrix}] 
  Let us compute the connection matrix $E_1$. Using the representations
  \eqref{Phi_0} and \eqref{Phi_1} for $\Phi^{(0)}(z)$ and $\Phi^{(1)}(z)$ we
  have
\begin{multline}
\label{proof:E_1}
\Phi_{11}^{(1)}(z)\,E_1^{11}+\Phi_{12}^{(1)}(z)\,E_1^{21}\\
\shoveleft{=\frac{\beta+\tau}{\delta}\,z^{\alpha}\,(z-1)^{\tau}\,
\F(a,b;a+b-c+1;1-z)\,E_1^{11}
}\\
\shoveright{+z^{\alpha}\,(z-1)^{-\tau}\,\F(c-a,c-b;c-a-b+1;1-z)\,E_1^{21}
}\\
\shoveleft{=\,\Phi_{11}^{(0)}(z)\,
=\,-\frac{2\alpha+1}{\delta}z^{\alpha}\,(z-1)^{\tau}\,\F(a,b;c;z)
}\\
\shoveleft{\overset{(\ref{connect0_1})}{=}\,
-\frac{2\alpha+1}{\delta}z^{\alpha}\,(z-1)^{\tau}\,
\left[
\frac{\Gamma(c)\,\Gamma(c-a-b)}{\Gamma(c-a)\,\Gamma(c-b)}\,
\F(a,b;a+b-c+1;1-z)
\right.
}\\
\left.
\shoveright{+(1-z)^{c-a-b}\,\frac{\Gamma(c)\,\Gamma(a+b-c)}{\Gamma(a)\,\Gamma(b)}\,
\F(c-a,c-b;c-a-b+1;1-z)
\right]
}\\
\shoveleft{=\,-\frac{2\alpha+1}{\delta}\,
\frac{\Gamma(c)\,\Gamma(c-a-b)}{\Gamma(c-a)\,\Gamma(c-b)}\,
z^{\alpha}\,(z-1)^{\tau}\,\F(a,b;a+b-c+1;1-z)
}\\
-\frac{2\alpha+1}{\delta}\,
\frac{\Gamma(c)\,\Gamma(a+b-c)}{\Gamma(a)\,\Gamma(b)}\,
z^{\alpha}(z-1)^{\tau}(1-z)^{-2\tau}\,\F(c-a,c-b;c-a-b+1;1-z)\,,
\end{multline}
since Theorem~\ref{th_Fuchs} implies $c-a-b=-2\,\tau$.

Here the branches of $z$ and $(z-1)$ are fixed by
$0<\arg z<2\pi$ and $0<\arg (z-1)<2\pi$. Using
\begin{equation*}
(z-1)^{\tau}(1-z)^{-2\tau}=(z-1)^{-\tau} e^{2\tau\pi i}
\end{equation*}
and identifying the coefficients at the hypergeometric functions
in
\eqref{proof:E_1} we obtain
\begin{align*}
E_1^{11}\,&=\,-\frac{2\alpha+1}{\beta+\tau}\,
\frac{\Gamma(c)\Gamma(c-a-b)}{\Gamma(c-a)\Gamma(c-b)},\\
E_1^{21}\,&=\,-\frac{2\alpha+1}{\delta}\,
\frac{\Gamma(c)\Gamma(a+b-c)}{\Gamma(a)\Gamma(b)}e^{2\tau\pi i}.
\end{align*}
It is easy to verify that the equality
\begin{equation*}
\Phi_{21}^{(0)}(z)\,=\,\Phi_{21}^{(1)}(z)\,E_1^{11}
+\Phi_{22}^{(1)}(z)\,E_1^{21}
\end{equation*}
also holds. Similarly we can obtain the formulas for
$E_1^{12}$ and $E_1^{22}$.
The computation for $E_\infty$ is analogous.
\end{proof}


\end{document}

%% file: param.pstex_t
\begin{picture}(0,0)%
\epsfig{file=param.pstex}%
\end{picture}%
\setlength{\unitlength}{4144sp}%
\begingroup\makeatletter\ifx\SetFigFont\undefined%
\gdef\SetFigFont#1#2#3#4#5{%
  \reset@font\fontsize{#1}{#2pt}%
  \fontfamily{#3}\fontseries{#4}\fontshape{#5}%
  \selectfont}%
\fi\endgroup%
\begin{picture}(5120,2580)(439,-2629)
\put(3001,-751){\makebox(0,0)[lb]{\smash{\SetFigFont{9}{10.8}{\rmdefault}{\mddefault}{\updefault}\special{ps: gsave 0 0 0 setrgbcolor}$\tfrac{1}{2}+i\tfrac{\sqrt{3}}{2}$\special{ps: grestore}}}}
\put(3001,-2001){\makebox(0,0)[lb]{\smash{\SetFigFont{9}{10.8}{\rmdefault}{\mddefault}{\updefault}\special{ps: gsave 0 0 0 setrgbcolor}$\tfrac{1}{2}-i\tfrac{\sqrt{3}}{2}$\special{ps: grestore}}}}
\end{picture}